\documentclass[12pt]{amsart}
\usepackage{amsmath,amscd,amssymb,amsthm,array}

\usepackage{mathtools}

\usepackage{hyperref}

\usepackage{comments}
\newComments\DL{D}{red}
\newComments\Bj{Bj}{red}
\newComments\IS{IS}{magenta}

\newcommand {\Lra} {{\Longrightarrow}}
\newcommand {\Llra} {{\Longleftrightarrow}}

\usepackage{DLdef1}

\let\ssec\subsection

\renewcommand {\ssbegin}[2][*]
 {\refstepcounter{subsection}%
\if#1*
\addcontentsline{toc}{subsection}{\thesubsection.\hskip 1pc #2}%
\else
\addcontentsline{toc}{subsection}{\thesubsection.\hskip 1pc #2. #1}%
\fi
 \def \secno {\gdef \secno {}{\ssecfont
\thesubsection.\hskip 2ex}%
 }%
 \begin{#2}}

\renewcommand {\sssbegin}[2][*]
 {\refstepcounter{subsubsection}
\if#1*
\addcontentsline{toc}{subsubsection}{\thesubsubsection.\hskip 1pc #2}%
\else
\addcontentsline{toc}{subsubsection}{\thesubsubsection.\hskip 1pc #2. #1}
\fi
 \def \secno {\gdef \secno {}{\ssecfont \thesubsubsection.\hskip 2ex}%
 }%
 \begin{#2}}

\renewcommand {\parbegin}[2][*]
 {\refstepcounter{paragraph}
\if#1*
\addcontentsline{toc}{paragraph}{\theparagraph.\hskip 1pc #2}%
\else
\addcontentsline{toc}{paragraph}{\theparagraph.\hskip 1pc #2. #1}
\fi
 \def \secno {\gdef \secno {}{\ssecfont \theparagraph.\hskip 2ex}%
 }%
 \begin{#2}}

\newcommand {\wht}{{\text{wht}}}

\begin{document}

\title{Analogs of Bol  operators on superstrings}

\author{Sofiane Bouarroudj${}^{a, *}$, Dimitry Leites${}^{a,b}$, Irina~Shchepochkina${}^c$}
\address{${}^a$Division of Science and Mathematics\\
New York University Abu Dhabi\\
Po Box 129188, United Arab Emirates; sofiane.bouarroudj@nyu.edu
\\
${}^b$Department of Mathematics\\
Stockholm University, Stockholm, Sweden; dl146@nyu.edu, mleites@math.su.se\\
${}^c$Independent University of Moscow,
B. Vlasievsky per., d. 11, RU-119 002 Moscow, Russia; irina@mccme.ru\\
$^{*}$ The corresponding author}

\begin{abstract} The Bol operators are unary differential operators between spaces of weighted densities on the 1-dimensional manifold  invariant under projective transformations of the manifold. On the $1|n$-dimensional supermanifold (superstring) $\mathcal{M}$,
we classify 
analogs of Bol operators invariant under   
the simple maximal subalgebra $\mathfrak{h}$ of the same rank as its simple ambient superalgebra $\mathfrak{g}$ of vector fields on $\mathcal{M}$ and containing all elements of negative degree of $\mathfrak{g}$ in a~ $\mathbb{Z}$-grading. We also consider the Lie superalgebras of vector fields $\mathfrak{g}$ preserving a~ contact structure on the superstring~ $\mathcal{M}$. 
We have discovered many new operators. \end{abstract}



\maketitle

\markboth{\itshape Sofiane Bouarroudj\textup{,} Dimitry Leites\textup{,} Irina Shchepochkina}
{{\itshape Analogs of Bol  operators on superstrings}} 

\thispagestyle{empty}

\tableofcontents
\setcounter{tocdepth}{3}

\section{Introduction} \label{Intro}

\subsection{Veblen's problem} 
On a~given $m$-dimensional manifold $\cM$, let $T_U(V)$ or just $T(V)$ be the space of $V$-valued tensor fields over a~ domain $U\subset \cM$.
Let $J_A$ be the Jacobi matrix of the diffeomorphism $A$ computed in coordinates at points $X\in \cM$ and $A^{-1}(X)\in \cM$. Let $\rho:\GL(m)\tto\GL(V)$ be a~representation.  The group $\Diff \cM$ of diffeomorphisms of $\cM$ acts on $T(V)$ by invertible changes of coordinates as follows:
\[
A(t)(X):=\rho(J_A)(t(A^{-1}(X))\text{~~for any $A\in\Diff \cM$, $X\in \cM$, $t\in T(V)$}. 
\]

In 1928, at the IMC, O.~Veblen
formulated a~problem over $\Ree$ (see \cite{V}) which A.~Kirillov later reformulated in more
comprehensible terms of modern differential geometry as follows (see \cite{Kinv}):
\begin{equation}\label{Veb}
\begin{minipage}[l]{13cm}
``Describe operators $D: T(V_1)\otimes\cdots\otimes T(V_k)\longmapsto T(W)$\\  invariant under all invertible diffeomorphisms of $\cM$''. 
\end{minipage}
\end{equation}

For \textit{differential} operators, J.~Bernstein \cite{BL} further reformulated Veblen's problem  as a~ purely algebraic problem meaningful over any ground field $\Kee$. Thus interpreted, Rudakov's solution (see \cite{R1, R2}) of Veblen's problem is reviewed in \cite{GLS1} for several types of unary ($k=1$) operators: invariant under arbitrary diffeomorphisms, or under diffeomorphisms preserving either a~volume element, or a~symplectic form, or a~contact structure. 

Being purely algebraic and local (hence, $\cM$ is, actually, any open domain) Veblen's problem can be reformulated as follows and its generalization for  supermanifolds and passage to ground fields of characteristic $p>0$ are immediate. The answers, however, are not that easy to obtain, and there are many new operators to be discovered: for unary operators in the super case over $\Cee$, see \cite{GLS1}; for modular cases and $\dim\cM=1$, see \cite{BouL}. 

Observe that in presence of even indeterminates on a~ supermanifold on which functions with compact support or of fast decay at infinity are considered, the integral is an example of \textit{non-differential} invariant operator from the space of densities with compact support to the space of functions; for just one more known example of \textit{non-differential} invariant operator, first discovered by A.A.~Kirillov only by its symbol, see~ \cite{Kinv,IM}. 

\textbf{Hereafter the ground field is $\Cee$}. Let $\cF$ be the algebra of functions (which in our case are polynomials or formal series). Denote the Lie algebra of vector fields  by $\cL:=\fvect(m)=\fder\ \cF$. 

Let $\cL_0$ be the maximal subalgebra of finite codimension (for $\fvect(m|0)$ it consists of the vector fields that vanish at the origin). Define the Weisfeiler filtration 
\[
\cL=\cL_{-d}\supset\dots\cL_{-1}\supset\cL_0\supset \cL_1\supset \dots
\]
by letting  
$\cL_{-1}$ be a~minimal
$\cL_{0}$-invariant subspace strictly containing $\cL_{0}$; the other
terms being defined by the following formula, where $i\geq 1$ 
\begin{equation}
\label{1.3}
\cL_{-i-1}=[\cL_{-1}, \cL_{-i}]+\cL_{-i}\; \text{ and
}\; \cL_i =\{D\in \cL_{i-1}\mid [D,
\cL_{-1}]\subset\cL_{i-1}\}.
\end{equation}

Let $V$ be an $\cL_0$-module, such that $\cL_1V=0$, i.e., $V$ is, actually,
$\fgl(m)\simeq\cL_0/\cL_1$-module considered as an $\cL_0$-module. The space of \text{tensor fields of type} $V$ with formal coefficients is defined to be
\[
T(V):=\text{Hom}_{U(\cL_0)}(U(\cL), V)\simeq \Cee[[x]]\otimes V.
\]

To consider \textit{irreducible} $\cL_0/\cL_1$-modules $V_i$ and $W$, see \eqref{Veb}, is a~natural first step, especially in the super situation and over fields of positive characteristics. Obviously, it is natural to replace $\fg=\fvect(m)$ with any other simple Lie (super)algebra $\cL$ of vector fields with polynomial (or formal if we need closedness with respect to the $x$-adic topology) coefficients.

For the list of known such Lie superalgebras, see \cite{BGLLS}. It is also natural to consider operators~ \eqref{Veb} invariant under the action of a~maximal subalgebra $\fh$ of $\fg$.

 
\subsection{Feasible versions of Veblen's problem}\label{feas}  Let $k$ be the arity of the operator \eqref{Veb}. For $k=1$,  Veblen's problem was solved for several simple Lie superalgebras, see a~ review~ \cite{GLS1}.

For $k=2$, Grozman solved Veblen's problem for $\fvect(m)$, its divergence-free subalgebra $\fsvect(m)$, and partly for the Lie algebra of Hamiltonian vector fields, see \cite{Gr1, Gr2}. 

The super versions of Veblen's problem for $k=2$ are open for any Lie superalgebra $\cL$, apart from some occasional results for contact Lie superalgebras $\fk(2m+1|n)$ and $\fm(n)$ for small values of $k$ and $m|n$, see \cite{KLV, Lr}.

For $k\geq 3$, Veblen's problem is out of reach in general, though it is conceivable to find the operator $D$, see \eqref{Veb}, of any prescribed order for any $\cL$ and any arity $k$ for any given $V_i$ and $W$, where ``any'' is understood ``within reason" depending on computers available.

There are, however, cases where one can get the complete answer even for $k$ ``large'', as, e.g., in \cite{Bter, B06, Bc}. These are cases where the spaces $T(V_i)$ and $T(W)$ are \textit{weighted densities}, see \S\ref{WD}, i.e., rank 1 modules over the algebra $\cF$ of functions,  and moreover, the  operators $D$ are anti-symmetric, as in \cite{FF}.

There is also an exceptional case where the simple Lie superalgebra of vector fields $\cL$ is such that $\cL_0/\cL_1$ is solvable (and hence the irreducible $\cL_0/\cL_1$-modules are finite-dimensional). Over $\Cee$, this only happens if $\cL\simeq\fvect(1|1)$.

\subsection{The problems we consider} Here we consider Veblen's problem over $\Cee$ and solve it for finite-dimensional maximal simple $\Zee$-graded subalgebras $\fh$ of the vectorial Lie superalgebra $\fg=\fvect(m|n)$ or $\fk(2m+1|n)$ or $\fm(n)$ in the standard grading of $\fg$, and such that $\fh_-=\fg_-:=\oplus_{i<0}\fg_i$ and the maximal tori of $\fh$ and $\fg$ coincide.
\begin{equation}\label{notN}
\begin{minipage}[l]{14cm}
The adjective \textbf{simple} applied to $\fh$ and $\fg$ in this problem is not natural;  we should (and do) consider various ``relatives'' of simple algebras: their central extensions (not only nontrivial) and algebras of derivations.
\end{minipage}
\end{equation} 

In various problems of mathematical physics, it is natural to consider a~\textbf{maximal} simple subalgebra $\fs$ of a~simple vectorial Lie algebra $\fg$; for their classification, see \cite{LSh}, where ``simple''\ is understood in the above-mentioned ``weak'' sense, see \eqref{notN}. Let $i: \fs\longrightarrow \fg$ be the corresponding embedding.  The following lemma describes the scope of our examples; for notation used in this text, see \S~\ref{WD} and Subsection~\ref{elu}, for further details, see \cite{BGLLS}.

\sssbegin[The maximal simple subalgebras $\fs$ of the simple vectorial Lie superalgebra $\fg$ such that $\rk\fs=\rk\fg$, see \cite{LSh}]{Lemma}[\cite{LSh}]\label{MainL} The maximal simple subalgebras $\fs$ of the simple vectorial Lie superalgebra $\fg$ such that $\rk\fs=\rk\fg$ are as follows 
\begin{equation}\label{d=2}
\renewcommand{\arraystretch}{1.4}
\begin{tabular}{|l|}
\hline
$\fosp(n|2m+2)\subset \fk(2m+1|n)$ \\
\hline
$\fpgl(m|n+2)\simeq \fpgl(n+2|m)\subset \fk(2n+1|2m)$\\ 
\hline
$\fosp_a(4|2)\subset \fk(1|4)$\\
\hline$\fas\subset \fk\fas\subset \fk(1|6)$\\ 
\hline\hline\hline
$\fpgl(m+1|n)\simeq\fpgl(n|m+1)\subset\fvect(m|n)$ for $(m, n)\neq(0,1), (0,2)$\\
\hline\hline\hline
$\fpe(n+1)\subset\fm(n)$\\
\hline
$\fpgl(m+1|n+1)\simeq \fpgl(n+1|m+1)\subset \fm(n+m)$ for $m+n>0$ \\
\hline
$\fspe_{a,b}(n+1)\subset\fb_{a,b}(n)$ \\
\hline\end{tabular} \end{equation}
\end{Lemma}

Note that in Table \eqref{d=2} we always have \textit{projective} embeddings (modulo center), but since $\fpgl(a|b)\simeq\fsl(a|b)$ if $a\neq b$, we sometimes write $\fsl$ instead of $\fpgl$ for clarity of matrix realization.  

From Table \eqref{d=2} we see that there are three major types of embeddings: 

\textbf{type $\fk$}: Embeddings into Lie superalgebras of contact type $\fk$. 

\textbf{type $\fm$}:  Embeddings into Lie superalgebras of contact type $\fm(n)$ and its subalgebra $\fb_{a,b}(n)$;  according to Convention~\eqref{stringsOnly}, we consider only $n=1$ and do not consider the case of non-simple ambient algebra $\fb_{a,b}(1)$.

\textbf{type $\fvect$}: Embeddings  into Lie superalgebras of general type $\fvect$. They are considered in a~ separate paper.

\begin{equation}\label{stringsOnly}
\begin{minipage}[l]{13cm} 
Among embeddings into $\fk(2m+1|n)$, the case of $m=0$ is of particular interest: the supermanifold, on which we consider the problem, is what physicists call a~\textit{superstring}. Even in this case, we consider analogs of Bol operators only between spaces of weighted densities: as we show in \cite{BLS}, the case of operators between spaces of tensor fields on supermanifolds other than superstrings is out of reach, bar several exceptions.  \\ 
\end{minipage}
\end{equation}

We  leave the cases where $\rk\fs<\rk\fg$, also classified in \cite{LSh}, for the future together with the exceptional embedding ${\fas\subset \fk\fas\subset \fk(1|6)}$.

\sssbegin[Bol's theorem \cite{Bol}]{Theorem}[Bol's theorem \cite{Bol}]\label{Bol1} Let $m|n=1|0$, see line~  $1$ in eq.~\eqref{d=2}. The only nonzero and non-scalar $i(\mathfrak{sl}(2))$-invariant linear
differential operators are as follows.

For every $k\in\Nee:=\{1,2,\dots\}$, there exists a~unique $i(\mathfrak{sl}(2))$-invariant
order $k$ linear differential operator $B_k: {\cal F}_{(1-k)/2} \tto {\cal F}_{(1+k)/2}$; explicitly
\[
B_k(\varphi(dt)^{(1-k)/2})=\varphi^{(k)}(dt)^{(1+k)/2}.
\]
\end{Theorem}

These operators $B_k$ are called the \textit{Bol operators}; they were first described in \cite{Bol}. Every Bol operator determines a~ projective structure on $S^1$, i.e., a~ (local) identification of the circle with the projective line, see \cite{BO} and \cite[p.9]{OT}. A similar identification of a~ ``super circle'' with a~ superstring takes place in the super case, cf. \cite{GT,H,CMZ}, where a~ simplest super case is considered (there are 4 series of superstrings and several exceptions, see \cite{GLS} where the Lie superalgebras of the supergroups of diffeomorphisms of supercircles are described and classified).

Considering the spaces of polynomial weighted densities on 1-dimensional manifold as Verma modules with lowest weight vector of weight $\mu\in\Cee$, the description of Bol operators is a~ trivial corollary of the classification of irreducible $\fsl(2)$-modules (if $\mu\in\Cee\setminus \Zee_{\leq 0}$) and their finite-dimensional submodules (if $\mu\in \Zee_{\leq 0}$), see any classification of irreducible $\fsl(2)$-modules with lowest or highest weight vector in terms of Verma modules, e.g., \cite{Ber, Kir}.

\sssbegin[Bol's analogs for $i(\fosp (1|2))$-invariant operators]{Theorem}[Bol's analogs for $i(\fosp (1|2))$-invariant operators, \cite{G, GT}]\label{k11} Let $m|n=0|1$, see line~  $1$ in eq.~\eqref{d=2}. The only non-zero non-scalar  $i(\fosp(1|2))$-invariant differential operators ${\cal F}_{a}\tto{\cal F}_{b}$
are as follows.

For every $k\in\Zee_+:=\Nee\cup\{0\}$,  there is an $i(\fosp(1|2))$-invariant differential operator
\[
B_k: {\cal F}_{-k} \tto {\cal F}_{k+1} \text{~~given by~~} \varphi \alpha_1^{-\frac{k}{2}} \mapsto D_\theta \circ\partial_t^k(\varphi)\alpha_1^{\frac{k+1}{2}}.
\]
\end{Theorem}

\ssec{Super Bol operators $\cF_\lambda\tto \cF_\mu$ in terms of $\theta$'s and $\xi\eta$'s} In \cite{BaL}, the authors considered the realization of $\fk(1|n)$ preserving the distribution singled out by the contact form 
\be\label{BaLcf}
dt +\sum \theta_i d\theta_i.
\ee
Then,  see \eqref{KD},
\[
K_{\theta_i}= \theta_i \partial_t - \partial_{\theta_i}, \ \ D_{\theta_i}=\partial_{\theta_i} + \theta_i \partial_t.
\]

\sssbegin[Bol's analogs for $i(\fosp (n|2))$-invariant operators]{Theorem}[Bol's analogs for $i(\fosp (n|2))$-invariant operators, \cite{BaL}]\label{kBaL} The only non-zero non-scalar $i(\fosp(n|2))$-invariant differential operators $\cF_\lambda\tto \cF_\mu$ are the above-described for $n=0$ and $1$, and \textup{(see \eqref{KD})}
\[
B_k: \varphi\alpha_1^{\lambda/2}\mapsto D_{\theta_1}\circ\dots \circ D_{\theta_n}\partial_t^k(\varphi)\alpha_1^{\mu/2}\text{~~for any $k\in\Zee_+$},
\]
where  
\be\label{lamu}
\lambda=1-k-n 
\text{ and } \mu=k+1.
\ee
\end{Theorem}

Observe that in \cite{BaL}, the authors define $\cF_a$ not as we do \eqref{wden}, but as $\cF\vvol^a$, but denote $\vvol$ by $\alpha$, hence the values of $\lambda$ and  $\mu$  in \cite{BaL}   differ from the above. 


\sssbegin[A $\xi\eta$-version of \cite{BaL}]{Theorem}[A $\xi\eta$-version of \cite{BaL}, Eq. (2.1), see Th.~\ref{kBaL}]\label{ksieta} In the realization of $\fk(1|n)$ preserving the distribution given by the contact form \eqref{2.2.6}, for $n>2$, the only non-zero non-scalar $i(\fosp(n|2))$-invariant operators $B_k: \cF_\lambda\tto \cF_\mu$   
of parity $n\pmod 2$ are as follows for any $k\in\Zee_+$
\be\label{nOd}
\begin{array}{ll}
\varphi 
\mapsto D_\theta\circ [D_{\eta_1}, D_{\xi_1}]_+ \circ \cdots \circ[D_{\eta_{\tilde n}}, D_{\xi_{\tilde n}}]_+ \circ \partial_t^k(\varphi)
&\text{for $n=2 \tilde n+1$ odd},\\
\varphi 
\mapsto [D_{\eta_1}, D_{\xi_1}]_+ \circ \cdots \circ[D_{\eta_{\tilde n}}, D_{\xi_{\tilde n}}]_+ \circ \partial_t^k(\varphi)
&\text{for $n=2\tilde n$ even},
\end{array}
\ee
where \textup{(we use notation \eqref{comm} and \eqref{KD})}
\be\label{lamu1}
\lambda=1-k-n \text{ and } \mu=k+1.
\ee
\end{Theorem}


\ssbegin[New Bol's analogs for $i(\fosp (2|2))$-invariant operators]{Theorem}[New Bol's analogs for $i(\fosp (2|2))$-invariant operators]\label{k12} Let $m|n=0|2$, see line~  $1$ in eq.~\eqref{d=2}.

\emph{(i)} For every $k\in\Zee_+$, the only non-zero non-scalar even $i(\fosp(2|2))$-invariant differential operators are
\[
B_k: \varphi \alpha^{\lambda}\beta^{\mu}\mapsto \bar B(\varphi) \alpha^{\nu}\beta^{\omega}
\]
where
\be\label{bolsN2}
\lambda=-\nu= \nfrac{k+1}2, \ \ \mu=\omega=- c(k+1)  \text{~~for any $c\in\Cee$}
\ee
and where \textup{(we use notation \eqref{comm} and \eqref{KD})}
\[
\bar B_k=[D_\xi,D_\eta]_+ D_1^k+cD_1^{k+1} .
\]
\textup{(For $\mu=\omega=c=0$, we recover the result of \cite{BaL}.)}

\emph{(ii)} There are only two $1$-parameter families of non-zero non-scalar odd $i(\fosp(2|2))$-invariant operators \textup{(see \eqref{KD})} for any $\lambda\in\Cee$:  
\[
\begin{array}{ll}
B_{\xi, \lambda}: {\cal F}_{\lambda, -\lambda} \tto {\cal F}_{1+\lambda, -1-\lambda}, &
\bar B_{\xi, \lambda}(\varphi)=D_\xi(\varphi);\\[2mm]
B_{\eta, \lambda}: {\cal F}_{\lambda, \lambda} \tto {\cal F}_{1+\lambda, 1+\lambda}, &
\bar  B_{\eta, \lambda}(\varphi)=D_\eta(\varphi).
\end{array}
\]
\end{Theorem}

Note an extremely beautiful formula \eqref{?}.

\ssbegin[New Bol's analogs of $i(\fosp_a(4|2))$-invariant operators]{Theorem}[New Bol's analogs of $i(\fosp_a(4|2))$-invariant operators]\label{T4} Let $m|n=0|4$, see line~  $3$ in eq.~\eqref{d=2}.

\emph{(i)}  Let $a\neq -1, \infty$. The only non-zero non-scalar $i(\fosp_a(4|2))$-invariant differential operators are
\[
B_k: {\cal F}_{-(k+3)} \tto {\cal F}_{k+1} \text{~~given by~~}  
\varphi \alpha_1^{-(k+3)/{2}} \mapsto \bar B_k(\varphi) \alpha_1^{(k+1)/{2}}  \text{for any $k\geq 0$}, 
\]
where \textup{(we use notation \eqref{comm} and \eqref{KD})}
\be\label{?}
\bar B_k(\varphi)=  [D_{\eta_1}, D_{\xi_1}]_+ \circ [D_{\eta_2}, D_{\xi_2}]_+ \circ \partial_t^k(\varphi) +\nfrac{4(1-a)}{(k+2)(1+a)} \partial_t^{k+2}(\varphi).
 \ee

\emph{(ii)} The only non-zero non-scalar $i(\fosp_{-1}(4|2))$-invariant differential operators are 
\[
B_k: {\cal F}_{\lambda} \tto {\cal F}_{\lambda +2k} \text{~~given by~~} \bar B_k: \varphi \alpha_1^{\lambda/{2}} \mapsto \partial_t^k(\varphi)\alpha_1^{{\lambda}/{2} + k}
 \text{for any $k\geq 0$}.
\]

\emph{(iii)}  The only non-zero non-scalar  $i(\fosp_{\infty}(4|2))$-invariant differential operators are
\[
B_k: {\cal F}_{-(k+3)} \tto {\cal F}_{k+1} \text{~~given by~~}    
\varphi \alpha_1^{-(k+3)/{2}} \mapsto \bar B_k(\varphi) \alpha_1^{(k+1)/{2}}  \text{for any $k\geq 0$}, 
\]
where \textup{(the coefficient of $\partial_t^{k+2}$ is  $\mathop{\lim}\limits_{a\tto \infty}\nfrac{4(1-a)}{(k+2)(1+a)}$, see eq.~\eqref{?})}
\[
\bar B_k(\varphi)=  [D_{\eta_1}, D_{\xi_1}]_+ \circ [D_{\eta_2}, D_{\xi_2}]_+ \circ \partial_t^k(\varphi) -\nfrac{4}{(k+2)} \partial_t^{k+2}(\varphi).
 \]
\end{Theorem}

\ssec{Comments} In $\fk(1|0)$, there is only one simple maximal subalgebra: $\fsl(2)$; for the result, see Theorem \ref{Bol1}. 

In $\fk(1|1)$, there is only one simple maximal subalgebra: $\fosp(1|2)$; for the result, see  Theorem \ref{k11}. 

In $\fk(1|2)\simeq\fm(1)$, there is one simple maximal subalgebra: $\fosp(2|2)\simeq\fsl(1|2)\simeq\fsl(2|1)$ corresponding to the first 2 lines in Table~\eqref{d=2},  and
one maximal subalgebra ``close to simple'': $\fpe(2)\subset\fm(1)$ corresponding to line~  6  in Table~\eqref{d=2}; for the result, see Theorem \ref{bolPe}. 

In order to explain where this subalgebra $\fpe(2)$ comes from, 
recall ``occasional''\ isomorphisms, see the bottom line in eq.~\eqref{iso} where the Lie superalgebras are considered as \textbf{abstract}, whereas as filtered or $\Zee$-graded they are \textbf{not} isomorphic:
\begin{equation}\label{iso}
\begin{array}{l}
\fosp(2|2)\simeq\fsl(1|2)\simeq\fvect(0|2),\\
\fk(1|2)\simeq\fvect(1|1)\simeq\fm(1).\\
\end{array}
\end{equation}

Considering the Lie superalgebras in the second line of eq.~\eqref{iso} not as abstract, which is unreasonable in Veblen's problem, but as filtered, we see that the $(\cF, \fg)$-\textbf{bi}modules of rank~ 1 over $\cF$, where $\cF$ is the algebra of functions on $\cM$ and $\fg$ is a~Lie superalgebra of certain vector fields on $\cM$, are quite different:
over $\fk(1|2)$ and $\fm(1)$ we have 2-parameter bimodules $\cF_{a,b}$, whereas the modules of weighted densities over $\fvect(1|1)$ can only depend on 1 parameter, see \cite{BLS}.

Accordingly, in the cases $\fk(1|2)$, $\fvect(1|1)$, and $\fm(1)$, we respectively replace $\fsl(2)$ with the Lie superalgebra $\fh:=\fosp(2|2)$, $\fsl(1|2)$, and $\fpe(2)$, considered as $\Zee$-graded, and embedded in the $\Zee$-graded ambient algebra $\fg$ so that $\Zee$-grading is preserved and $i(\fh_-)=\fg_-$.
This forces us to consider $\fpe(2)$ as a~subalgebra of~ $\fm(1)$. 

\ssec{Dualization and permutation of arguments}\label{dual} For tensor fields with compact support, we have, as is easy to prove by integration,
\be\label{fakedual}
(T(W))^*= T(W^*)\otimes_\cF \Vol,\text{~~where $\Vol:=T(\tr)$} 
\ee 
and where  $\tr:A\tto \tr(A)$ for any $A\in \fgl(m|n)$ denotes both the (super)trace map and the $\fgl(m|n)$-module it defines. Let us define the dual spaces of tensor fields for the fields with polynomial or formal coefficients  via  formula \eqref{fakedual} as if the fields were with compact support. For any invariant operator ${D: T(V)\tto T(W)}$, the dual operator 
\[
D^*: (T(W))^*:=T(W^*)\otimes_\cF \Vol\tto (T(V))^*:=T(V^*)\otimes_\cF \Vol
\]
 is also invariant, and since $\Vol:=T(\tr)$, we can describe $D^*$ as 
\be\label{dual}
D^*: T(W^*\otimes \tr)\tto T(V^*\otimes \tr).
\ee
In particular, for the spaces of weighted densities, see Subsection~\ref{VWD}, let $D: \cF_a\tto \cF_b$ be an invariant operator, and $\Vol=\cF_v$, where $v$ depends on $\sdim \cM=m|n$. Then, ${D^*: \cF_{v-b}\tto \cF_{v-a}}$ must be (and is) invariant. This helps to double-check the answers, see \cite{Gr1, Gr2, BGL, BouL}.

\section{Weighted densities}\label{WD} The $\fvect(m|n)$-module $\Vol^a$ of weighted densities of weight $a\in\Kee$ is a~ rank~ $1$ module over the algebra of functions $\cF$ in indeterminates $x:=(u_1,\dots ,u_m;\xi_1, \dots ,\xi_n)$ generated by $\vvol^a$, where $\vvol$ is the volume element with constant coefficients, i.e., the class of the element ${du_1\wedge\dots \wedge du_m\partial_{\xi_1}\dots \partial_{\xi_n}}$ in the quotient of 
\[
E^m_\cF((\fvect(m|n))^*)\otimes_\cF S^n_\cF(\fvect(m|n))
\]
--- the tensor product over $\cF$ of the $m$th exterior power of the $\cF$-module of covector fields by the $n$th symmetric power of the $\cF$-module of vector fields --- modulo the maximal $\fvect(m|n)$-submodule of corank 1. The  $\fvect(m|n)$-action is given for any $D\in\fvect(m|n)$, $a\in\Kee$, and $f\in\cF$, by the Lie derivative, where  $p$  designates the parity function 
\be\label{Lvol}
L_D(f\vvol^a)=(D(f)+(-1)^{p(f)p(D)}af\Div(D))\vvol^a .
\ee

\ssec{Weighted densities on supermanifolds with a~ contact structure} There are two types of contact structures; in these cases, \textbf{certain} modules of tensor fields $T(V)$ on a~supermanifold $\cM$ form rank 1 modules over the algebra $\cF$, see \cite{BGLLS}. The corresponding Lie superalgebras are:

$\bullet$ $\fk(2n+1|m)$ preserving the distribution on $\cM$ given by an \textbf{odd} contact form $\alpha_1$ with \textbf{even} ``time'', see eq. \eqref{2.2.6}; 

$\bullet$ $\fm(r)$ preserving the distribution on $\cM$ given by an \textbf{even} contact form $\alpha_0$ with \textbf{odd} ``time'', see eq. \eqref{alp0}. 

\sssec{Case $\fk(2n+1|m)$}  It is more natural to express the $\fk(2n+1|m)$-modules of weighted densities  in terms of the contact form $\alpha_1$, see eq. \eqref{2.2.6}, than in terms of $\vvol$:
\be\label{wden}
\cF_\lambda:=\begin{cases}\cF\alpha_1^{\lambda/2}&\text{if $n+m>0$,}\\
\cF\alpha_1^{\lambda}&\text{if $n=m=0$}.\\
\end{cases}
\ee
with the $\fk(2n+1 |m)$-action given by the formula (compare with eq.~\eqref{Lvol} and \eqref{2.3.5}) 
\be\label{aNeVol}
L_{K_f}(\varphi \alpha_1^{\lambda})=(K_f(\varphi)+\lambda K_1(f) \varphi ) \alpha_1^\lambda.
\ee

\textbf{Important}: as $\fk(2n+1 |m)$-modules, $\Vol^a\simeq\cF_{a(2+2n-m)}$. 

$\bullet$ Let $u=(t, p_1, \dots , p_n, q_1, \dots , q_n)$ be even indeterminates,  
$\Theta=\begin{cases}(\xi, \eta)& \text{ if }\ m=2k\\(\xi, \eta,
\theta)&\text{ if }\ m=2k+1,\end{cases}$ be odd indeterminates
and let  
\begin{equation}
\label{2.2.6} \alpha_1=dt+\mathop{\sum}\limits_{1\leq i\leq
n}(p_idq_i-q_idp_i)+ \mathop{\sum}\limits_{1\leq j\leq
k}(\xi_jd\eta_j+\eta_jd\xi_j)+
\begin{cases}0&
\text{ if }\ m=2k\\
\theta d\theta&\text{ if }\ m=2k+1.\end{cases}  
\end{equation}

Consider the distribution singled out by the \textit{Pfaff equation} \index{Pfaff equation} 
\be\label{alpha1}
\text{$\alpha_1(X)=0$ for $X\in \fvect(2n+1|m)$.}
\ee
The Lie superalgebra, called \textit{contact} Lie superalgebra,  that preserves this distribution is
\begin{equation}
\label{k} \fk (2n+1|m)=\{ D\in \fvect (2n+1|m)\mid
L_D \alpha_1=f_D \alpha_1\text{ for some }f_D\in \Cee [t, p, q,
\theta]\}.
\end{equation}

\sssec{Case $\fm(n)$}   Similarly, set $u=q=(q_1, \dots , q_n)$,
let $\theta=(\xi_1, \dots , \xi_n; \tau)$ be odd. Set
\begin{equation}
\label{alp0}
\begin{array}{c}
\alpha_0=d\tau+\mathop{\sum}\limits_i(\xi_idq_i+q_id\xi_i),
\end{array}
\end{equation}
and call this form  \textit{pericontact}.

Consider the distribution singled out by the Pfaff equation
\be\label{alpha0}
\text{$\alpha_0(X)=0$ for $X\in \fvect(n|n+1)$.}
\ee
The Lie superalgebra, called \textit{pericontact Lie superalgebra},   preserves this distribution is
\begin{equation}
\label{m} \fm (n)=\{ D\in \fvect (n|n+1)\mid L_D\alpha_0=f_D\cdot
\alpha_0\text{ for some }\; f_D\in \Cee [q, \xi, \tau]\}.
\end{equation}
To call $\fm$ ``odd contact
superalgebra", as some authors do, is misleading since the even part of $\fm$ is nonzero.
To say ``the ``odd" contact superalgebra'' is OK, though a~bit long.

\subsubsection{Generating functions: contact series and their subalgebras}\label{SS:2.4.1} A laconic way to describe $\fk$,
$\fm$ and their subalgebras is via generating functions. \index{Function, generating} We express the contact fields of series $\fk$ (resp. $\fm$) in natural bases
$K_f$ (resp. $M_f$).

$\bullet$ \underline{Odd form $\alpha_1$}. For any $f\in\Cee [t, p,
q, \theta]$, set:\index{$K_f$, contact vector field} \index{$H_f$,
Hamiltonian vector field}
\begin{equation}
\label{2.3.1} K_f=(2-E)(f)\pder{t}-H_f + \pderf{f}{t} E,
\end{equation}
where $E=\sum_i y_i \pder{y_{i}}$ (here the $y_{i}$
are all the coordinates except $t$) is the \textit{Euler operator},\index{$E$, the Euler operator}
and $H_f$ is the hamiltonian field with Hamiltonian $f$ that
preserves $d\alpha_1$:
\begin{equation}
\label{2.3.2'} H_f=\mathop{\sum}\limits_{i\leq n}\left
(\pderf{f}{p_i} \pder{q_i}-\pderf{f}{q_i} \pder{p_i}\right)
-(-1)^{p(f)}\left(\pderf{f}{\theta} \pder{\theta}+\mathop{\sum}\limits_{j\leq
k}\left( \pderf{f}{\xi_j} \pder{\eta_j}+ \pderf{f}{\eta_j}
\pder{\xi_j}\right)\right).
\end{equation}

 $\bullet$ \underline{Even form $\alpha_0$}. For any $f\in\Cee [q,
\xi, \tau]$, set:
\begin{equation}
\label{2.3.3} M_f=(2-E)(f)\pder{\tau}- Le_f -(-1)^{p(f)}
\pderf{f}{\tau} E, 
\end{equation}
where $E=\sum_iy_i \pder{y_i}$ (here the $y_i$ are
all the coordinates except $\tau$), and\index{$E$, the Euler operator}
\begin{equation}
\label{2.3.4} Le_f=\mathop{\sum}\limits_{i\leq n}\left(
\pderf{f}{q_i}\ \pder{\xi_i}+(-1)^{p(f)} \pderf{f}{\xi_i}\
\pder{q_i}\right).
\end{equation}
\index{$M_f$, pericontact vector field} \index{$Le_f$, periplectic
vector field} Since
\begin{equation}
\label{2.3.5}
\renewcommand{\arraystretch}{1.4}
\begin{array}{l}
 L_{K_f}(\alpha_1)=2 \pderf{f}{t}\alpha_1=K_1(f)\alpha_1, \\
L_{M_f}(\alpha_0)=-(-1)^{p(f)}2 \pderf{
f}{\tau}\alpha_0=-(-1)^{p(f)}M_1(f)\alpha_0,
\end{array}
\end{equation}
it follows that $K_f\in \fk (2n+1|m)$ and $M_f\in \fm (n)$. (These vector fields actually span the respective Lie superalgebras.)
Observe that
\[
p(Le_f)=p(M_f)=p(f)+\od.
\]

$\bullet$ To the (super)commutators $[K_f, K_g]$ or $[M_f, M_g]$
there correspond \textit{contact brackets}\index{Poisson
bracket}\index{Contact bracket}\index{$\{-, -\}_{k.b.}$, contact bracket}\index{$\{-, -\}_{m.b.}$, pericontact bracket} of the generating functions:
\[
\renewcommand{\arraystretch}{1.4}
\begin{array}{l}
{}[K_f, K_g]=K_{\{f, \; g\}_{k.b.}};\qquad
 {}[M_f, M_g]=M_{\{f, \;
g\}_{m.b.}}.\end{array}
\]
The explicit formulas for the contact brackets are as follows. Let
us first define the brackets on functions that do not depend on $t$
(resp. $\tau$).

The \textit{Poisson bracket} $\{-,-\}_{P.b.}$ 
in the realization with the form $\omega_0$ 
is given by the formula
\begin{equation}
\label{pb}
\renewcommand{\arraystretch}{1.4}
\begin{array}{ll}
 \{f, g\}_{P.b.}&=\mathop{\sum}\limits_{i\leq n}\ \bigg(\displaystyle\pderf{f}{p_i}\
\pderf{g}{q_i}-\ \pderf{f}{q_i}\
\pderf{g}{p_i}\bigg)\\
&-(-1)^{p(f)}\bigg(\displaystyle\pderf{f}{\theta}\ \pderf{
g}{\theta}+
\mathop{\sum}\limits_{j\leq m}\left(\displaystyle\pderf{f}{\xi_j}\ \pderf{ g}{\eta_j}+\pderf{f}{\eta_j}\
\pderf{ g}{\xi_j}\right)\bigg)\\
&\text{ for any }f, g\in \Cee [p, q, \xi, \eta,
\theta].
\end{array}
\end{equation}

The \textit{Buttin bracket} $\{-,-\}_{B.b.}$ \index{Antibracket $=$ Buttin
bracket $=$ Schouten bracket} \index{Buttin
bracket $=$ Schouten bracket  $=$ Antibracket}\index{Schouten bracket$=$ Antibracket $=$ Buttin
bracket} is given by the formula
\begin{equation}
\label{2.3.7} \{ f, g\}_{B.b.}=\mathop{\sum}\limits_{i\leq n}\
\bigg(\pderf{f}{q_i}\ \pderf{g}{\xi_i}+(-1)^{p(f)}\
\pderf{f}{\xi_i}\ \pderf{g}{q_i}\bigg)\text{ for any }f, g\in \Cee
[q, \xi].
\end{equation}

In terms of the Poisson and Buttin brackets,
respectively, the contact brackets are
\begin{equation}
\label{2.3.8} \{ f, g\}_{k.b.}=(2-E) (f)\pderf{g}{t}-\pderf{f}
{t}(2-E) (g)-\{ f, g\}_{P.b.}
\end{equation}
and
\begin{equation}
\label{2.3.9} \{ f, g\}_{m.b.}=(2-E)
(f)\pderf{g}{\tau}+(-1)^{p(f)} \pderf{f}{\tau}(2-E) (g)-\{ f,
g\}_{B.b.}. 
\end{equation}

\ssec{Divergences}

\sssec{Divergences on $\fk(2n+1 |m)$} For  $m\neq 2n+2$, the restriction of the divergence from $\fvect(2n+1 |n)$ to $\fk(2n+1 |m)$ is equal to
\[
\Div(K_f)= (2n+2-m)\partial_t(f).
\]
Therefore, $\fk(2n+1 |2n+2)\subset\fsvect(2n+1 |2n+2)$. If $m\neq 2n+2$, then the divergence-free subalgebra of
$\fk(2n+1 |m)$ is the \textit{Poisson} Lie superalgebra of functions independent of $t$, spanned by $K_f$ such that $\partial_t(f)=0$.

The \textit{standard grading} in $\cF$, where $\deg t=2$ and the degree of any other indeterminate being equal to 1, induces the \textit{standard grading} in $\fk(2n+1 |m)$: we set $\deg(K_f)=\deg(f)-2$.

Since codimension of $(\fk(2n+1 |m)_0)'$ in $\fk(2n+1 |m)_0$ is equal to 1, there is one (up to a~ non-zero factor) trace on $\fk(2n+1 |m)_0$. Hence, there should be one (up to a~ non-zero factor) divergence on $\fk(2n+1 |m)$, see \cite{BGLLS}. 
Set $\tr(K_t\vert_{\fg_{-1}})=1$, then the divergence is given by $\partial_t$ for any $n$ and $m$. 

\paragraph{The case $\fk(1 |2)$} This case is exceptional. The  algebra $\fk(1 |2)_0$ of degree 0 elements is commutative and $\dim\fk(1 |2)_0=2$, so there are 2 traces on $\fk(1 |2)_0$, and hence 2 divergences on $\fk(1 |2)$. 
Let $\alpha_1:=dt+\xi d\eta+\xi d\eta$. To the operators
\[
\tilde K_\xi(f)=(-1)^{p(f)}(\partial_\eta-\xi\partial_t)(f), \ \tilde K_\eta(f)=(-1)^{p(f)}(\partial_\xi-\eta\partial_t)(f)\text{~~and~~}\tilde K_1=\partial_t
\]
there corresponds the second divergence
\begin{equation}\label{Div_2}
\Div_2:=\tilde K_\eta\tilde K_\xi-\tilde K_1.
\end{equation}

Let $\beta$ be the symbol of the class of a~pseudo-differential form that can be ``visualized'' as $\beta:=\sqrt{[d\xi(d\eta)^{-1}]}$. Its actual expression is irrelevant, we need to know only how the Lie derivative acts on it.  Let $\wht$ be the weight with respect to $K_t$ and $K_{\xi\eta}$; then 
\be\label{wht}
\text{$\wht(\xi)=(1,1)$ and $\wht(\eta)=(1,-1)$, so $\wht([d\xi(d\eta)^{-1}])=(0,2)$.}
\ee
 Since it is natural to have a~ unit of weight be equal to 1, not 2, we take the square root.

The Lie derivative acts as follows 
\begin{equation}\label{L_Dk12}
\begin{array}{l}
L_{K_f}(\alpha_1^a\beta^b)=(2a\partial_t(f)+ (-1)^{p(f)} b \Div_2(K_f))(\alpha_1^a\beta^b)\\
=((2a- (-1)^{p(f)} b )\partial_t(f)+ (-1)^{p(f)} b\tilde K_\eta\tilde K_\xi(f))(\alpha_1^a\beta^b).
\end{array}
\end{equation}

The space $\cF_{a,b}$ of weighted densities over $\fk(1|2)$ is a~rank 1 module over the algebra of functions generated by $\alpha^{a/2}_1\beta^{b}$. (The above description is OK unless the characteristic of the ground field is equal to $2$, cf.~\cite{BGLLS}.)

\sssec{Case $\fm(n)$} 
On $\fm(n)$, there are two analogs of the divergence on $\fvect$. One is the operator $\Div$ --- the restriction of the divergence from $\fvect(n|n+1)$ to $\fm(n)$:
\[
L_{M_f}\vvol=\Div(M_f)\vvol;
\]
the other one, $\Div_\tau:=2\partial_\tau
\circ \text{sign}$, arises from the Lie derivative acting on $\alpha_0$:
\begin{equation}\label{divonm}
L_{M_f}\alpha_0 =2(-1)^{p(f)}\partial_\tau(f)\alpha_0.
\end{equation}
Explicitly, (here $E=\sum y_i\partial_{y_i}$ the sum of all the coordinates except $\tau$)
\begin{equation}\label{DeltaM_f1}
\Div(M_f)=2(-1)^{p(f)}\Delta^\fm(f), \text{where $\Delta^\fm:=(1-E)\partial_\tau-\sum_i\frac{\partial^2}{\partial_{\xi_i}\partial_{q_i}}$}.
\end{equation}

Rank 1 modules $\cF_{a,b}$, where $a,b\in\Cee$, over the algebra $\cF$ are generated by $\alpha_0^a\vvol^b$, on which $\fm(n)$ acts as follows: 
\begin{equation}\label{DeltaM_f2}
L_{M_f}(\alpha_0^a\vvol^b) =2(-1)^{p(f)}(a\partial_\tau(f)+ b~\Delta^\fm(f))(\alpha_0^a\vvol^b).
\end{equation} 

\ssec{The volume elements and  weighted densities}\label{VWD} In this subsection all isomorphisms are up to the change of parity functor $\Pi$ because (a) there are two natural choices of parity on the space $\Vol$, see \cite{MaG}, and (b) neither choice affects the invariance.

Over $\fvect(m|n)$, we have $\Vol:=\cF\vvol$. Since $\vvol=[dx_1\dots dx_m\partial_{\xi_1}\dots \partial_{\xi_n}]$, then $\Vol\simeq\cF_{m-n}$, where the weight is taken relative the Euler operator $E=\sum x_i\partial_{x_i}+\sum \xi_j\partial_{\xi_j}$. 

Over $\fk(2m+1|n)$, we have $\Vol\simeq\cF_{2m+2-n}$, where the weight is taken relative $K_t$. 

Over $\fm(n)$, we have $\Vol\simeq\cF_{-2}$, where the weight is taken relative $M_{\tau}$.

\section{The $\fosp(n|2m+2)\subset\fk(2m+1|n)$ cases} \label{BolsKosp} 

According to \eqref{stringsOnly}, in this paper we consider only the case $m=0$. 

\ssec{Preliminaries}  Recall the definition of $K_{f}$ and introduce the operators $D_{f}$, see \cite{Shch}, for $f=1, p, q$ (even) and $\xi, \eta, \theta$ (odd):  
\be\label{KD}
\begin{array}{l}
\text{$D_{1}=\partial_t$,\ \ $K_1=2\partial_t$};\\
\text{$D_{\eta_i}:=\xi_i \partial_t+\partial_{\eta_i}$,\ $D_{\xi_i}:=\eta_i \partial_t+\partial_{\xi_i}$, \ $D_{\theta}:=\theta\partial_t+\partial_{\theta}$;  \ $D_{p_j}= p_j \partial_t+\partial_{q_j}$, \ $D_{q_j}= q_j \partial_t-\partial_{p_j}$;}\\[2mm]
\text{$K_{\xi_i}=\xi_i \partial_t - \partial_{\eta_i} $,\ $K_{\eta_i}= \eta_i \partial_t-\partial_{\xi_i}$, \ $K_{\theta}= \theta\partial_t - \partial_{\theta}$; \ $K_{p_j}= p_j \partial_t-\partial_{q_j}$, \ $K_{q_j}= q_j \partial_t+\partial_{p_j}$}.
\end{array}
\ee

The cases $\fk(1|n)$ for $n=2$ and $n=4$ are exceptional: if $n=2$, the weight of densities acquires an extra parameter and new Bol-type operators associated with this parameter; if $n=4$, the embedded algebra can be deformed and there are new Bol-type operators depending on the parameter of deformation.

The invariance of the operator $B: \cF_a\tto \cF_b$ or, equivalently, of the map $B(f\alpha^a)=\bar B(f)\alpha^b$, where $\alpha$ denotes the generator of the $\cF$-module of weighted densities (i.e., $\alpha_1$, or $\alpha_0\vvol$, or $\vvol$) over $\fvect(\cM)$ or any its subalgebra $\fg$ means that
\be\label{invA1}
L_D(B(f\alpha^a))=(-1)^{p(B)p(D)}B(L_D(f\alpha^a)) \text{~~for any $D\in \fg$ and $f\in\cF$}.
\ee

Now, let $\fh$ be a~subalgebra of $\fg:=\fvect(\cM)$, or $\fk(\cM)$ or $\fm(\cM)$ on the supermanifold $\cM$ such that $\fh_-=\fg_-$ relative selected $\Zee$-gradings of $\fh$ and $\fg$. Then, 
\[
L_D(\alpha)=0 \text{~~for any $D\in \fg_-$},
\]
and therefore if condition \eqref{invA1} is satisfied, then the action of the operators $\bar B$ and $L_D$ on the space of functions supercommute. 

For the elements of $\fm(n)$, recall the definition of $M_{f}$ and introduce the operators $D_{f}$ for $f=1, q$ (even) and $\xi$ (odd): 
\be\label{MD}
\begin{array}{l}
\text{$D_{1}=\partial_\tau$,\ \ $M_1=2\partial_\tau$};\\
\text{$M_{\xi_i}:=\xi_i \partial_\tau+\partial_{q_i}$,\ $M_{q_i}:=q_i \partial_\tau-\partial_{\xi_i}$;}\\[2mm]
\text{$D_{\xi_i}:=\xi_i \partial_\tau-\partial_{q_i}$,\ $D_{q_i}:=q_i \partial_\tau+\partial_{\xi_i}$}.
\end{array}
\ee

Recall that the \textit{super commutator} and the \textit{super anti-commutator} of operators $X$ and $Y$ are, respectively,
\be\label{comm}
\text{$[X, Y]:=XY+YX$ and $[X, Y]_+:=XY-YX$.}
\ee

\ssbegin[All invariant differential operators have constant coefficients in an invariant  basis \cite{Shch}]{Theorem}[All invariant differential operators have constant coefficients in an invariant   basis \cite{Shch}]\label{Peetre}
 Let $\fh$ be a~ subalgebra  in a~ vectorial Lie superalgebra $\fg$, and $\fh_-=\fg_-$. 
 
 Let $\fg:=\fk(2m+1|n)$. Any $\fh_-$-invariant differential operator $\bar B$ is a~polynomial with constant coefficients in the operators~$D_f$, see eq. \eqref{KD}.  
 
If $\fg=\fm(n)$, then any $\fh_-$-invariant differential operator $\bar B$ is a~polynomial with constant coefficients in the operators $D_f$, see eq. \eqref{MD}.

If $\fg=\fvect(m|n)$, then any $\fh_-$-invariant differential operator $\bar B$ is a~polynomial with constant coefficients in partial derivatives.
\end{Theorem}

\begin{proof}  Let $\fg:=\fk(2m+1|n)$. Let us arrange monomials in $D_{\xi_i}$, $D_{\eta_i}$, and $D_1=\partial_t$ in this order to form a~ basis  in the $\cF$-module of 
differential operators. (Since $[D_{\xi_i}, D_{\eta_i}]=\pm D_1$, any basis contains $D_1$ and one of the two: either $D_{\xi_i}D_{\eta_i}$ or $D_{\eta_i}D_{\xi_i}$.) All operators $D_{f}$, where $\deg f<2$, commute with~ $\fg_-$. Then,  for any monomial $F$ in the operators $D_{f}$, any $\varphi\in\cF$ and $K_g\in\fg_-$, we have
\[
[K_g, \varphi F]=K_g(\varphi)F.
\]
Hence, any $\fg_-$-invariant differential operator must have constant coefficients ($K_g(\varphi)=0$).

\begin{equation}\label{invComp}
\begin{minipage}[l]{14cm} 
Since the condition \eqref{invA1} should hold for any function $f\in\cF$, then having selected a~ basis of $\cF$ homogeneous with respect to a~ grading, we see that all components of $\bar B$ homogeneous with respect to the same grading should be invariant.
\end{minipage}
\end{equation}

Thanks to \eqref{invComp}, we can assume that $\bar B$ is homogeneous of degree $s$. Since 
\[
K_t(f)=(\deg f)f, \text{~and~} K_t(\alpha_1^a)=2a\alpha_1^a,
\]
then \eqref{invA1} for $D=K_t$ turns into 
\[
\begin{array}{ll}
&BK_t(f\alpha_1^a)=(\deg f+2a)\bar B(f)\alpha_1^b\\
=&K_tB(f\alpha_1^a)=(\deg \bar B(f)+2b)\bar B(f)\alpha_1^b=(s+\deg f+2b)\bar B(f)\alpha_1^b,
\end{array}
\]
i.e., 
\[
2a=s+2b.
\]

Analogously, 
\be\label{invV}
L_D(\vvol)=(\Div D)\vvol.
\ee

Over $\fk(2m+1|n)$, we define the spaces of weighted densities generated, as modules over the algebra of functions $\cF$, by 
$\alpha_1^a$ (or $\alpha_1^a\beta^b$ in dimension $1|2$), see \eqref{aNeVol}, NOT by $\vvol^a$, as in \cite{BaL}. 

Observe also that $L_{K_{\xi_i\eta_i}}(\alpha_1)=0$ for all $i$, and hence the actions of the operators $\bar B$ and $K_{\xi_i\eta_i}$ on the space of functions must commute. From commutation relations
\[
\begin{array}{l}
{}[K_{\xi_i\eta_i}, D_{\xi_i}]=D_{\xi_i}, \ \ [K_{\xi_i\eta_i}, D_{\eta_i}]=-D_{\eta_i},\\
{}[K_{\xi_i\eta_i}, D_{\xi_j}]=[K_{\xi_i\eta_i}, D_{\eta_j}]=0\text{~~for $i\neq j$}
\end{array}
\]
it follows that $\bar B$ should contain products $D_{\xi_i}D_{\eta_i}$ and $D_{\eta_i}D_{\xi_i}$. 

Thus, $\bar B$ should be a~ sum of monomials in $[D_{\xi_i},D_{\eta_i}]_+$ and $\partial_t$ with constant coefficients. 

\textbf{This argument is not applicable to the $\fk(1|2)$-invariant operators $\cF_{a,b}\tto\cF_{c,d}$ with $b, d\neq 0$. Odd Bols can and do appear, see Theorem~\ref{k12}}. Actually, there are odd operators $\cF_{a,b}\tto\cF_{c,d}$ invariant with respect to the whole $\fk(1|2)$, not only Bols, see \cite{Lr}.

For $\fg=\fm(n)$ and $\fg=\fvect(m|n)$, the proof is analogous.
\end{proof}

\sssec{Recapitulation from \cite{Lr}: elucidation of Theorem \ref{k12}}\label{sssRec} For the a~ description of $\fg=\fk(1|2)$-invariant differential operators $D: T(V)\tto T(W)$ with irreducible $\fg_0$-modules $V$ and $W$, see \cite{Lr}.  In this case, all modules $T(V)$ are of the form $\cF_{a,b}$, so the answers in  \cite{Lr} should be among our Bols. Let us reproduce the results from \cite{Lr} and amend them to get Bols. 

On the $n$-dimensional manifold $M$, we have $(\Omega^i)^*=\Omega^{n-i}$. Recall that on the supermanifold $\cM$, the dual to the space of \textit{differential forms} $\Omega^i:=T(E^i(V))$, where $E^i$ is the operator of raising to the $i$th exterior power and $V$ is the tautological $\fgl(V)=\fgl(\dim(\cM))$-module, is the space of \textit{integrable forms}
\[
\Sigma_{-i}(\cM):=(\Omega^i(\cM))^*\otimes_\cF \Vol(\cM)=T(E^i(V^*)\otimes \tr), 
\]
where $\Vol(\cM)$ is the $\cF$-module of volume forms (aka Berezinian). There are two complexes (series of mappings such that the product of two consecutive mappings is equal to 0)
\be\label{d}
\begin{array}{l}
\Omega^0 \stackrel{d}{\tto}\Omega^{1} \stackrel{d}{\tto}\dots \Omega^i \stackrel{d}{\tto}\Omega^{i+1}\stackrel{d}{\tto}\dots ; \\
\dots\stackrel{d^*}{\tto} \Sigma_j \stackrel{d^*}{\tto} \Sigma_{j+1}\stackrel{d^*}{\tto}\dots
\Sigma_{-1} \stackrel{d^*}{\tto} \Sigma_{0}.
\end{array}
\ee

In \cite{Lr}, following Rudakov, instead of modules $T(V):=\Hom_{U(\fg_0\oplus\fg_+)}(U(\fg), V)$, Leites considered the induced modules 
\[
I(V):=U(\fg)\otimes_{U(\fg_0\oplus\fg_+)}V\simeq\Cee[\fg_-]\otimes V=(T(V^*))^*,
\]
where $\fg_-:=\oplus_{i<0}\fg_i$ and $\fg_+:=\oplus_{i>0}\fg_i$, and $V$ is considered as a~ $\fg_0\oplus\fg_+$-module such that $\fg_+V=0$. Since $\fg_0$ is commutative, then any irreducible $\fg_0$-module is 1-dimensional. Observe that
\be\label{relK12}
(K_\xi)^2=\nfrac12[K_\xi, K_\xi]=0,\ \ (K_\eta)^2=\nfrac12[K_\eta, K_\eta]=0,\ \ [K_\xi, K_\eta]=K_1.
\ee
The invariant operators $D: T(V)\tto T(W)$ are in one-to-one correspondence with  invariant operators 
\[
D^*: T(W)^*=I(W^*)\tto T(V)^*=I(V^*), 
\]
and hence \textit{singular vectors} $f\in I(V^*)$, i..e.,  such that $\deg f<0$ and $\fg_+f=0$.
Since $\fg_+$ is generated by 3 operators $K_{t\xi}$, $K_{t\eta}$, and $K_{t\xi\eta}$, it suffices to solve systems of 3 equations for every $\deg f\leq -2$ and only 2 equations if $\deg f=-1$ since in this case $K_{t\xi\eta}f=0$ holds automatically.

In what follows, let  $K_tv=Ev$ and $K_{\xi\eta}v=Hv$ for some $E, H\in \Cee$ and $v\in V$. 

$\bullet$ $\deg f=-1$. Then,  $f=(aK_\xi  +bK_\eta)v$, where $a,b\in\Cee$, and $ab=0$ because $f$ is supposed to be a~ weight vector with respect to $K_t$ and $K_{\xi\eta}$. Recall that $\fg_+v=0$. Let us solve the system
\be\label{1}
\begin{array}{ll}
(A)&K_{t\xi}f=0;\\
(B)&K_{t\eta}f=0.
\end{array}
\ee
The system \eqref{1} is equivalent to
\be\label{11}
\begin{array}{ll}
(A)&b(K_{\xi\eta}-K_t)v=0, \\ 
(B)&a(K_{\xi\eta}+K_t)v=0.
\end{array}
\ee
Therefore, recall \eqref{wht},
\be\label{calc}
\begin{array}{ll}
(A)&\text{$H=E=\lambda$ and $\wht(f)=(\lambda-1,\lambda-1)$},\\
(B)&\text{$H=-E=\lambda$ and $\wht(f)=(\lambda-1,-\lambda+1)$}.
\end{array}
\ee

Let us try to understand the nature of these operators. We know that the exterior differential $d: \cF_{0,0}\tto \Omega^1$ is invariant with respect to any changes of variables, not only $\fg$-invariant. Since $\cF\alpha_1$ is a~ $\fg$-invariant submodule in $\Omega^1$, we have projections 
\[
\begin{array}{l}
B_{\xi, 0}: =\pr_\xi\circ d: \cF_{0,0}\tto \cF[d\xi]=\cF_{1,1};\\
B_{\eta, 0}: =\pr_\eta\circ d: \cF_{0,0}\tto \cF[d\eta]=\cF_{1,-1},
\end{array}
\]
where $[d\xi]$ and $[d\eta]$ are classes of $d\xi$ and $d\eta$, respectively, modulo $\cF\alpha_1$.
Calculations \eqref{1} and \eqref{calc} show  that there are $\fk(1|2)$-invariant analogs of these operators
\[
\text{$B_{\xi, a}: =\cF_{a,a}\tto \cF_{a+1,a+1}$ and 
$B_{\eta, a}: =\pr_\eta\circ d: \cF_{a,a}\tto \cF_{a+1,a-1}$.}
\]

Since $i(\fosp(2|2))$ contains $K_{t\xi}$ and  $K_{t\eta}$, the operators $B_{\xi, a}$ and 
$B_{\eta, a}$ are Bols.

Several  $\fk(1|2)$-invariant unary differential operators of the degree 2 found in \cite{Lr} are obvious: consider the restrictions of the composition
\be\label{-1}
\Sigma_{-1}\stackrel{d}{\to}\Sigma_{0}\stackrel{\phantom{aa}\vvol^{-1}}{\simeq}\Omega_{0}\stackrel{d}{\to}
\Omega_{1}
\ee
acting from submodules of the form $\cF_{a,b}$ to the quotients of the same form. 

There are many more Bol operators than restrictions of operators \eqref{-1} to rank 1 submodules of weighted densities, see \eqref{bolsN2}. Unlike operators \eqref{-1}, we can not interpret them.

$\bullet$ $\deg f=-3$. Then,  $f=(aK_\xi K_1  +bK_\eta K_1)v$, where $a,b\in\Cee$, and $ab=0$. It is easy to check that $f$ can not be annihilated by both $K_{t\xi}$ and $K_{t\eta}$. The same is true (by induction) for any vector of degree $-(2n+1)<-1$.

\sssec{Proof of Theorem \ref{k12}}\label{sssPf} (i)
Clearly, any \textbf{even} operator commuting with $\fg_-$ is of the form
\begin{equation}\label{operA}
B:\vf\alpha^\lambda\beta^\mu\tto   \bar B(\vf)\alpha^\nu\beta^\omega, \text{ where  $\bar B=aD_\xi D_\eta\partial_t^k+b\partial_t^{k+1}$ and $a,b\in\Cee$}.
\end{equation}

The operator $K_f$ commutes with the operator (\ref{operA}) if and only if  for any  $\vf$ we have
\begin{equation}\label{comf}
K_f\bar B(\vf)+(2\partial_t(f)\nu+\omega\Div_2(f))\bar B(\vf)=\bar B K_f(\vf)+2\lambda \bar B(\partial_t(f)\vf)+\mu \bar B(\Div_2(f)\vf) .
\end{equation}

To find out when the operator (\ref{operA}) commutes with the subalgebra $\fosp(2|2)\subset\fk(1|2)$, it suffices to find when the operator \eqref{operA} commutes with $K_{t^2}$ and $\fosp(2|2)_-=\fk(1|2)_-$ because they generate $\fosp(2|2)$. The corresponding calculations are not that easy (especially by bare hands), so let us first consider commutations with $K_t$ and $K_{\xi\eta}$. This will help in further calculations.

Set
\[
\gr (t)=2, \quad \gr (\xi)=\gr(\eta)=1,
\]
let $\deg(x)=1$ for any $x=t, \xi$ or $\eta$.

The action of $K_t$:
\[
\begin{array}{l}
K_t=2t\partial_t+\xi\partial_\xi+\eta\partial_\eta, \; K_t\alpha=2\alpha,\; K_t \beta=0, \; K_t(\vf)=\gr(\vf)\vf,\\
K_t(\bar B\vf)=(\gr(\bar B)+\gr(\vf))\bar B\vf=(-2k-2+\gr(\vf))\bar B\vf.
\end{array}
\]
Thus, the equation (\ref{comf}) for $f=t$ is equivalent to the equation
\begin{equation}\label{lanu}
\gr(\bar B)+\gr(\vf)+2\nu=\gr(\vf)+2\lambda\ \Llra \  2\nu+2k+2=2\lambda\ \Llra \  \nu=\lambda+k+1.
\end{equation}

Now, consider the action of $K_{\xi\eta}=-\xi\partial_\xi+\eta\partial_\eta$:
\[
\begin{array}{l}
K_{\xi\eta}\alpha=0,\; K_{\xi\eta} \beta=\beta,\\
{}[K_{\xi\eta},D_\xi]=D_\xi,\; [K_{\xi\eta},D_\eta]=-D_\eta, \; [K_{\xi\eta},D_1]=0 \Lra [K_{\xi\eta},\bar B]=0.
\end{array}
\]

Hence, for $f=\xi\eta$ the equation (\ref{comf}) is of the form
\begin{equation}\label{muom}
\omega \bar B(\vf)=\mu \bar B(\vf) \ \Llra \  \omega=\mu.
\end{equation}

Now, pass to $K_{t^2}$ and introduce $T$:
\begin{equation}\label{Kt2}
2T:=K_{t^2}=2t^2\partial_t+2t(\xi\partial_\xi+\eta\partial_\eta), \; K_{t^2}\alpha=4t\alpha,\; K_{t^2}\beta=2\xi\eta\beta,\; T(\vf)=(\deg\vf ) t\vf.
\end{equation}

Accordingly, the equations (\ref{comf}) for $f=t^2$ with equality (\ref{muom})  taken into account take the form
\[
T\bar B(\vf)+(2t\nu+\omega\xi\eta)\bar B(\vf)=\bar B T(\vf)+2\lambda \bar B(t\vf)+\omega \bar B(\xi\eta\vf),
\]
which taking  (\ref{Kt2}) into account, is equivalent to
\begin{equation}\label{comT}
T\bar B(\vf)+(2t\nu+\omega\xi\eta)\bar B(\vf)=(\deg\vf+2\lambda) \bar B(t\vf)+\omega \bar B(\xi\eta\vf).
\end{equation}


Consider the action of $\bar B$ on monomials:
\begin{equation}\label{barA1}
\bar B(t^s)= 
    \begin{cases}
    0 &\text{for $s\le k$}\\
    (a+b)(k+1)! &\text{for $ s=k+1$}\\
    \frac{s!}{(s-k-1)!}\left( (a+b)t^{s-k-1}-a(s-k-1)\xi\eta t^{s-k-2}\right) &\text{for $ s>k+1$}\\
    \end{cases}
\end{equation}

\begin{equation}\label{barA2}
 \begin{array}{l}
\bar B(\xi t^s)= \begin{cases}
    0 &\text{for $ s\le k$}\\
    b\frac{s!}{(s-k-1)!}\xi t^{s-k-1} &\text{for $ s>k$}\\
    \end{cases}\\
  \bar B(\eta t^s)=\begin{cases} 0 &\text{for $  s\le k$}\\
                                    (2a+b)\frac{s!}{(s-k-1)!}\eta t^{s-k-1} &\text{for $  s>k$}\\
\end{cases}\end{array}
\end{equation}

\begin{equation}\label{barA3}
\bar B(\xi\eta t^s)=\begin{cases}
    0 & \text{for $  s< k$}\\
    -ak! &\text{for $   s=k$}\\
    \frac{s!}{(s-k)!}\left(-at^{s-k}+(a+b)(s-k)\xi\eta t^{s-k-1}\right) &\text{for $   s>k$}\\
    \end{cases}
\end{equation}

Let us write the condition (\ref{comT}) for $\vf=t^k$ taking (\ref{barA1}) into account. We get
\[
(k+2\lambda)(a+b)(k+1)!+\omega(-ak!)=0,
\]
implying the condition for $\omega$:
\begin{equation}\label{om}
a\omega=(k+1)(a+b)(k+2\lambda)
\end{equation}

Now, let us write the condition  (\ref{comT}) for  $\vf=t^{k+1}$ taking (\ref{barA1}) into account:
\[
\begin{array}{l}
(2t\nu+\omega\xi\eta)(a+b)(k+1)!\\
=(k+1+2\lambda)(k+2)!\left( (a+b)t-a\xi\eta\right)+\omega (k+1)!(-at+(a+b)\xi\eta)
\end{array}
\]
After simplification we get
\[
t\left((a+b)(2\nu-(k+2)(k+1+2\lambda))+\omega a\right) +a(k+1+2\lambda)(k+2)\xi\eta =0,
\]
which is equivalent to the system of two conditions
\[
a\omega=(a+b)((k+2)(k+1+2\lambda)-2\nu) \text{ and } a(k+1+2\lambda)=0.
\]

Substituting $\nu=\lambda+k+1$ in these conditions and taking (\ref{lanu}) into account we see that the first condition is equivalent to (\ref{om}), while the second one is equivalent to the following condition
\[
\text{$a =0$ or $ \lambda =\nu$.}
\]

If $a\ne 0$, we can set $a=1$. Then,  the condition $ \lambda =\nu$ with (\ref{lanu}), (\ref{muom}) and (\ref{om}) taken into account yields the following values of parameters:
\begin{equation}\label{par}
\lambda=-\nfrac{k+1}2, \quad \nu= \nfrac{k+1}2, \quad \mu=\omega=-(k+1)(1+b).
\end{equation}

Actually, this is the answer desired. To be absolutely sure, we have to

1) check that for $a=0$ the operator $B$ is not invariant;

2) check that for $a=1$ if conditions (\ref{par}) are satisfied, then $B$ is indeed invariant, i.e., conditions (\ref{comT}) hold for any  function $\vf$.

Let $a=0$.  Then, we can assume that $\bar B=\partial_t^k$. In this case, equalities $\xi\eta\bar B(\vf)=\bar B(\xi\eta\vf)$  and  (\ref{comT}) turn into
\begin{equation}\label{a0}
T\bar B(\vf)+2\nu t\bar B(\vf)=(\deg \vf+2\lambda)\bar B(t\vf).
\end{equation}

Let $\vf=\vf_0 t^k$, where the function $\vf_0$ depends only on $\xi, \eta$ and $\deg \vf_0=d$.
Then, the equation (\ref{a0}) takes the form 
\[
\begin{array}{l}
k!(d+2\nu)\vf_0t=(d+k+2\lambda)(k+1)!\vf_0t\\
 \ \Llra \  d+2\nu=(d+k+2\lambda)(k+1)\ \Llra \  dk=2\nu-(k+2\lambda)(k+1).
\end{array}
\]
This equality should hold for any $d=0,1,2$, which is impossible since RHS is a~ constant. Thus, the case $a=0$ is ruled out.

It only remains now to verify the invariance of $B$ for parameters given by relations (\ref{par}). Just plug (\ref{par}) in (\ref{comT}) and verify for all monomials  $\varphi$ using formulas (\ref{barA1})--(\ref{barA3}). 

Thus for $\vf=t^s$ and $s>k+1$ we see that the LHS of (\ref{comT}) is of the form 
\[
\begin{array}{l}
\frac{s!}{(s-k-1)!}((1+b)(s-k-1)t^{s-k}-(s-k-1)(s-k)\xi\eta t^{s-k-1})\\
+t(k+1)\frac{s!}{(s-k-1)!}((1+b)t^{s-k-1}-(s-k-1)\xi\eta t^{s-k-2})-\\
-(k+1)(1+b)\xi\eta \frac{s!}{(s-k-1)!}(1+b)t^{s-k-1}\\
=\frac{s\cdot s!}{(s-k-1)!}(1+b)t^{s-k}-\left(\frac{(s+1)!}{(s-k-2)!}+\frac{s!(k+1)(1+b)^2}{(s-k-1)!}  \right)\xi\eta t^{s-k-1}
\end{array}
\]
whereas the RHS of (\ref{comT}) is of the form
\[
\begin{array}{l}
(s-k-1)\frac{(s+1)!}{(s-k)!}((1+b)t^{s-k}-(s-k)\xi\eta t^{s-k-1})-
\\
-(k+1)(1+b)\frac{s!}{(s-k)!}(-t^{s-k}+(1+b)(s-k)\xi\eta t^{s-k-1})=
\\
=(1+b)\frac{s!}{(s-k)!}((s-k-1)(s+1)+k+1)t^{s-k}-\left(\frac{(s+1)!}{(s-k-2)!}+\frac{s!(k+1)(1+b)^2}{(s-k-1)!}  \right)\xi\eta t^{s-k-1}.
\end{array}
\]

Since $(s-k-1)(s+1)+k+1=s(s-k)$, then LHS=RHS. Thus, on monomials $\vf=t^s$ the equality (\ref{comT}) holds.

Now, check the equality (\ref{comT}) for $\vf=\xi t^s$. In this case, $\xi\eta\bar B(\vf)=\bar B(\xi\eta\vf)=0$. The LHS of (\ref{comT}) is of the form
\[
(s-k+k+1)b\nfrac{s!}{(s-k-1)!}\xi t^{s-k}=\nfrac{(s+1)!}{(s-k-1)!}b\xi t^{s-k},
\]
whereas the RHS is of the form
\[
(s+1-k-1)b\nfrac{(s+1)!}{(s-k)!}\xi t^{s-k}=\nfrac{(s+1)!}{(s-k-1)!}b\xi t^{s-k}.
\]
Thus, the equality (\ref{comT}) holds for such monomials as well.

Calculations for $\vf=\eta t^s$ are analogous.

It remains to verify equality (\ref{comT}) for $\vf=\xi\eta t^s$. In this case, $\bar B(\xi\eta\vf)=0$, but $\xi\eta\bar B(\vf)\ne 0$. The RHS is equal to
\[
\begin{array}{l}
(s+2-k-1)\frac{(s+1)!}{(s-k+1)!}(-t^{s-k+1}+(s-k+1)(1+b)\xi\eta t^{s-k})\\
=\frac{(s+1)!}{(s-k)!}(-t^{s-k+1}+(s-k+1)(1+b)\xi\eta t^{s-k}).
\end{array}
\]
Let us calculate the LHS:
\[
\begin{array}{l}
\frac{s!}{(s-k)!}(-(s-k)t^{s-k+1}+(s-k+1)(s-k)(1+b)\xi\eta t^{s-k})
\\
 +(k+1)\frac{s!}{(s-k)!}(-t^{s-k+1}+(s-k)(1+b)\xi\eta t^{s-k}) - \omega\frac{s!}{(s-k)!}\xi\eta t^{s-k}
\\
=-(s+1)\frac{s!}{(s-k)!}t^{s-k+1}+\left(\frac{s!}{(s-k-1)!} (1+b)(s+2)+\frac{s!}{(s-k)!}(1+b)(k+1) \right)\xi\eta t^{s-k}
\\
=-\frac{(s+1)!}{(s-k)!}t^{s-k+1}+\frac{s!}{(s-k)!}(1+b)((s+2)(s-k)+k+1)\xi\eta t^{s-k}.
\end{array}
\]
LHS=RHS since 
\[
(s+2)(s-k)+k+1=(s+1)(s-k+1).
\] 

Thus, the equality (\ref{comT}) is verified for all monomials $\vf$, and hence the operator 
\[
B: \vf\alpha^{\lambda}\beta^{\mu}\longmapsto  \bar B(\vf)\alpha^{\nu}\beta^{\omega},
\]
is invariant, where
\[
\bar B=D_\xi D_\eta\partial_t^k+b\partial_t^{k+1}, \quad \nu=-\lambda=\nfrac{k+1}2, \quad \mu=\omega=-(k+1)(1+b).
\]

Set $c:=1+b$. Then
\[
\bar B=(D_\xi D_\eta-\partial_t)\partial_t^k+c\partial_t^{k+1}.
\]
Since $[D_\xi,D_\eta]=2\partial_t$, we have
\[
-D_\eta D_\xi=D_\xi D_\eta-2\partial_t~\Lra~ [D_\xi,D_\eta]_+=2(D_\xi D_\eta-\partial_t).
\]
Since $B$ is defined up to a~ non-zero scalar factor, having denoted $2B$ by $B$ and $2c$ by $c$, we get the final form of $B$:
\[
B: \vf\alpha^{-\nfrac{k+1}2}\beta^{-c(k+1)}\longmapsto ([D_\xi,D_\eta]_+\partial_t^k+c\partial_t^{k+1})(\vf)\alpha^{\nfrac{k+1}2}\beta^{-c(k+1)}.
\]

 (ii) \textbf{Odd operators}. For the classification of $\fk(1|2)$-invariant differential operators of degree~1, see \cite{Lr} and the proof reproduced in Subsection \ref{sssRec}. There are no $\fk(1|2)$-invariant differential operators of order $>2$, see \cite{Lr}. There are no Bols either, because already invariance with respect to $K_{t\xi}$ and $K_{t\eta}$
is violated. \qed

\sssec{An explanation of the results of Theorem \ref{T4}}\label{elu} Recall that  Kaplansky discovered a~parametric family of deformations of $\fosp(4|2)$. He denoted the members of the family  by $\Gamma(\sigma_1,\sigma_2,\sigma_3)$, see \cite{Kapp}. Here, we use a~ short and suggestive notation $\fosp_a(4|2)$ or  $\fosp(4|2;\alpha)$, where the parameter $a$ was introduced in formulas \eqref{basis} whereas $\alpha$ is described as follows. 

For $\fg=\fosp(4|2;\alpha)$,  we have $\fg_\ev=\fsl_1(2)\oplus\fsl_2(2)\oplus\fsl_3(2)$ and $\fg_\od:=\id_1\otimes \id_2\otimes \id_3$, where the isomorphic copies of $\fsl(2)$ and their respective tautological modules are numbered for convenience. The Lie superalgebra structure on $\fg$ is completely defined by the bracket of any two odd elements, i.e., by the mapping $S^2(\fg_\od)\longrightarrow \fg_\ev$ defined up to a~permutation of the isomorphic summands $\fsl(2)$ of $\fg_\ev$ 
 or, which is the same, on non-vanishing simultaneously coefficients $\sigma_i$, where $i=1,2,3$, with which the projections $S^2(\fg_\od)\tto\fsl_i(2)$ enter the bracket. The Jacobi identity is satisfied if and only if $\sigma_1+ \sigma_2+ \sigma_3=0$, see \cite{Kapp}. Since the $\sigma_i$ are defined up to a~ common non-zero multiple and at least two of them must be non-zero, let $\sigma_1=1$ and $\alpha:=\nfrac{\sigma_i}{\sigma_j}$, where $\sigma_j\neq 0$ for $\{i,j\}=\{2,3\}$.
Thus, $S_3$ acts on the space of parameter $\alpha$. For explicit conditions when $\fosp(4|2;\alpha_1)\simeq \fosp(4|2;\alpha_2)$, see \cite{BGL}.  The isomorphisms $\fosp(4|2;\alpha)\simeq \fosp(4|2;\alpha')$ are generated by the
transformations:
\begin{equation}\label{osp42symm}
\alpha\longmapsto -1-\alpha\ , \qquad \alpha\longmapsto
\nfrac{1}{\alpha}\ ,
\end{equation}
so the other isomorphisms are 
\begin{equation}\label{osp42symm1}
\alpha\longmapsto
-\nfrac{1+\alpha}{\alpha}\ ,\quad\alpha\longmapsto
-\nfrac{1}{\alpha+1}\ ,\quad\alpha\longmapsto -\nfrac{\alpha}{\alpha+1}.
\end{equation}
Unlike
$\fosp(4|2;\alpha)$ or $\fosp_a(4|2)$, the notation 
$D(2,1;\alpha)$, though often used, is ill-chosen as explained, e.g., in \cite{CCLL}: the root system $D$ is not distinguished among the three possibilities.

We know, see \cite{BGL}, that $S_3$ acts on the parameter ~$\alpha$ of $\fosp(4|2;\alpha)$ with 6 elements in the orbit of general position and 3 elements in the exceptional cases. The reader might wonder if we lost some of the operators corresponding to the other values of parameter of the isomorphic algebras. Let us show that we did not lose anything. For brevity, we say ``$f$ acts'' meaning ``$K_f$ acts'', etc.

\paragraph{On $i(\fosp_a(4\vert 2))\subset\fk(1\vert 4)$} In Section~\ref{elu} we explain the difference between the parametric family $\fosp(4|2;\alpha)$ considered in \cite{BGL} and the family of embedded algebras $\fosp_a(4|2)$ considered here. 

The following generating functions form a~ basis of $i(\fosp_a(4|2))\subset\fk(1|4)$, see \cite{LSh}: 
\be\label{basis}
\begin{array}{ll}
\fosp_a(4|2)_{\leq 0}&=\fk(1|4)_{\leq 0};\\
\fosp_a(4|2)_{1}:& (1+a)t \xi_i +(a-1)\partial_{\eta_i} \Xi \text{ and } (1+a)t \eta_i +(a-1)\partial_{\xi_i} \Xi, \text{ where } i=1,2;\\
\fosp_a(4|2)_{2}:&\frac{1}{2}(1+a)t^2-(a-1)\Xi, \text{ where } \Xi=\xi_1 \xi_2 \eta_1 \eta_2.\end{array}
\ee

Observe that $i(\fosp_1(4|2))\simeq\fosp(4|2;1)\simeq\fosp(4|2)$, so the list of Bol operators in this case should be as in Theorem~\ref{ksieta}. It is, but $\fosp(4|2;1)\simeq\fosp(4|2; -2)\simeq\fosp(4|2; -\nfrac12)$, see \cite{BGL}, whereas the answers for $a=-2$ and $-\nfrac12$ are different, see Theorem \ref{T4}. Is there a~ mistake? 

Not easy to see immediately why, but there is no mistake: the corresponding spaces of parameters $\alpha$ and $a$ have non-isomorphic symmetry groups.

Observe also that $\fosp(4|2;\alpha)$ is not simple for $\alpha =0, -1$, and the list of Bols might differ in these cases. But we are working with $\fosp_a(4|2)$, not with $\alpha$, and for it the Bols are described by the general formula \eqref{?} evaluated at $a=0$ and as $a\to \infty$ as well, which we gave separately for clarity. The algebra $\fosp_{-1}(4|2)$ is not simple and the Bols are predictably different.

The case of $a=-1$ is exceptional: as spaces, $i(\fosp_{-1}(4|2))\simeq{\Lambda(\xi_1,\xi_2,\eta_1,\eta_2)\oplus \Cee t}$. 

In $\fg$, the subalgebra generated by $\{f\in\Lambda(\xi, \eta)\mid \int f\vvol(\xi, \eta)=0\}$  
is an ideal  isomorphic to $\fpo'(0|4)$; the elements $t$ and $\Xi$, see eq.~\eqref{basis}, act on the ideal by outer derivations. The Lie superalgebra $\fpo'(0|4)$ has an ideal generated by $\Cee\cdot 1$, the quotient is $\fh'(0|4)\simeq \fpsl(2|2)$. The ``additional'' to $\fsl(2|2)$ elements $1, t, \Xi$ form a~subalgebra NOT isomorphic to $\fsl(2)$, because $\{1,\Xi\}=0$.
 
If $a\ne -1$, set
\[
b:=\nfrac{a-1}{a+1}, \quad A_i:=t\xi_i+b\nfrac{\partial\Xi}{\partial\eta_i}, \; B_i:=t\eta_i+b\nfrac{\partial\Xi}{\partial\xi_i} \;(i=1,2), \quad F:=-t^2+2b\Xi.
\]

Then, the non-zero brackets of $[\fg_{-1},\fg_{-1}]$ are only $\{\xi_i,\eta_i\}=1$, and the non-zero brackets of $[\fg_{1},\fg_{1}]$ are only $\{A_i,B_i\}=F$. Besides,
\[
\{1,A_i\}=2\xi_i, \; \{1,B_i\}=2\eta_i \quad \text{ and } \{F,\xi_i\}=2A_i, \; \{F,\eta_i\}=2B_i.
\]
Let us now look how $\fg_0=\fo(4)\oplus \Cee t\simeq \fsl(2)_1\oplus\fsl(2)_2\oplus \Cee t$ acts. The element $t$ is a~grading operator $\ad_t|_{\fg_k}=k\cdot\id$.

The elements 
\[
\text{$x=u(\xi_1\eta_1-\xi_2\eta_2)+v \xi_1\eta_2+w\xi_2\eta_1$, where $u,v,w\in\Cee$}
\] span an ideal in $\fg_0$ isomorphic to $\fsl(2)$. We denote it $\fsl(2)_1$ for definiteness. The matrix of the restriction of $\ad_x$ to $\fg_{-1}$ in the basis $\xi_1, \xi_2,\eta_1,\eta_2$, and to $\fg_1$ in the basis $A_1,A_2, B_1, B_2$ is
\[
\begin{pmatrix}
-U^t & 0\\
0 & U\\
\end{pmatrix}, \text{ where~~} U=\begin{pmatrix}
        u & w\\
        v & -u\\
        \end{pmatrix} .
 \]
 
The elements
\[
\text{$y=\lambda(\xi_1\eta_1+\xi_2\eta_2)+\mu \xi_1\xi_2+\nu \eta_1\eta_2$, where $\lambda, \mu,\nu\in\Cee$}
\]
 also span an ideal in $\fg_0$ isomorphic to $\fsl(2)$. We denote it $\fsl(2)_2$. 
 
 The matrix of restriction of the operator $\ad_y$ to $\fg_{-1}$ in the basis $\xi_1, \xi_2,\eta_1,\eta_2$, and to $\fg_1$ in the basis $A_1,A_2, B_1, B_2$ is 
 \[
\begin{pmatrix}
  - \lambda & 0& 0 & -\mu \\
  0 & -\lambda & \mu & 0 \\
  0 &- \nu & \lambda & 0\\
  \nu & 0 & 0 & \lambda\\
  \end{pmatrix}.
\]
We see that the parameter $b$ arises only in commutators in $[\fg_{-1},\fg_1]$. Explicitly, they are (for $i,j=1, 2$):
\[
\begin{array}{l}
\{\xi_i,A_i\}=\{\eta_i,B_i\}=0,\\
\{\xi_i,A_j\}=(1+b)\xi_i\xi_j, \ \ \{\xi_i,B_j\}=(1-b)\xi_i\eta_j, \\ 
\{\eta_i,B_j\}=(1+b)\eta_i\eta_j, \ \ \{\eta_i,A_j\}=(1-b)\eta_i\xi_j \text{~~for $i\ne j$}.
\end{array}
\]
Finally,
\[
\begin{array}{ll}
\{\xi_1,B_1\}= -t+\xi_1\eta_1+b \xi_2\eta_2, & \{\xi_2,B_2\}=-t+b \xi_1\eta_1+\xi_2\eta_2,\\
\{\eta_1,A_1\}=-t-\xi_1\eta_1-b \xi_2\eta_2, & \{\eta_2,A_2\}=-t-b \xi_1\eta_1-\xi_2\eta_2. \end{array}
\]

We see that for $b=-1$, which is equivalent to $a=0$, we have $[\fg_{-1},\fg_1]\simeq\fsl(2)_1\oplus \Cee t$. So, $\fosp_0(4|2)$ contains an ideal $\fg_{-2}\oplus \fg_{-1}\oplus [\fg_{-1},\fg_1]\oplus \fg_1\oplus \fg_2$, isomorphic to $\fpsl(2|2)$, on which $\fsl(2)_2$ acts by outer derivations.

For $b=1$, which is equivalent to $a=\infty$, we have $[\fg_{-1},\fg_1]=\fsl(2)_2\oplus \Cee t$. So, $\fosp_{\infty}(4|2)$ contains an ideal $\fg_{-2}\oplus \fg_{-1}\oplus [\fg_{-1},\fg_1]\oplus \fg_1\oplus \fg_2$, isomorphic to a~copy of $\fpsl(2|2)$ on which $\fsl(2)_1$ acts by outer derivations.

Observe again 
that the values of parameter $a$ do not coincide with the same values of $\alpha$ from \cite{BGL}, except occasionally, e.g. for $a=\alpha=1$.

Furthermore, the group acting on the space of parameter $a$ and establishing an isomorphism $\fosp_{a_1}(4|2)\simeq\fosp_{a_2}(4|2)$ is not isomorphic to the symmetry group of the parameter $\alpha$. Indeed, having embedded $\fosp_{a}(4|2)$ in $\fk(1|4)$ in the case where $a \ne -1$, we distinguish one $\fsl(2)$, spanned by $1, t$ and $\frac{1}{2}(1+a)t^2-(a-1)\Xi$, see eq.~\eqref{basis}. This embedding endows $\fg_\od$ with a~$\Zee$-grading which shows how the brackets of two odd elements with values in this distinguished $\fsl(2)$ should look. 

The parameter $a$ arises \textbf{only} in commuting $\fg_{-1}$ with $\fg_1$ with values in $\fg_0$ which is now isomorphic to $\fsl(2) \oplus \fsl(2) \oplus \Cee$. So now, we can only permute two copies of $\fsl(2)$, and hence the group establishing isomorphisms is now $S_2$, not $S_3$. 

The formulas for commutators show that $S_2$ is generated by the mapping $b\mapsto -b$, or equivalently, $a\mapsto \frac 1a$. The values $a=\pm 1$ are exceptional points of this action: their orbits consist of one point, not of two.

Explicitly, the isomorphism is given by the map
\[
f: \vf(\xi_1,\xi_2,\eta_1,\eta_2)\mapsto \vf(\xi_1,\eta_2,\eta_1,\xi_2).
\]

\ssbegin{Remarks} 1) In \cite{BaL}, Basdouri and Laraiedh described Bol's operator in the case of $\fk(1|n)$ in realization preserving the distribution singled out by the form $dt+\sum \theta_id\theta_i$. This realization has advantage as compared with our realization in terms of $\xi$ and $\eta$: the formulas in \cite{BaL} is uniform for all~ $n$. Basdouri and Laraiedh used this uniformity
to give an analytic proof of their result, which is also an advantage as compared with computer-assisted proof which is justly considered less solid than theoretical. 

On the other hand, the $\xi\eta$ realization of $\fk(1|n)$ we considered makes it possible to use weight decomposition and facilitates computations (either by pencil or by computer).

2) Basdouri and Laraiedh did not consider the second parameter of the weight densities when $n=2$, and therefore missed two odd Bol operators because the main object of interest in \cite{BaL} was  
$H^1(\fk(1|n), \fosp(n|2); \cal F_\lambda)$, not  analogs of Bol operators. 

3) The ground field in \cite{BaL} is $\Ree$, but the answer is true for any field of characteristic zero.
\end{Remarks}

\section{The $\fpgl(m|n+2)\simeq \fpgl(n+2|m)\subset \fk(2n+1|2m)$ cases}\label{BolsKsl}

Let $x$ be indeterminates $m$ of which are even and $n$ are odd, and two more indeterminates $\xi$ and $\eta$ of the same parity. On $\fpgl(m|n+2)$ realized by linear differential operators, consider the grading given by
\[
\deg x_i=1 \text{~~for all~~}i;\ \ 
\deg \xi =0;\ \ 
\deg \eta =2.
\]

By convention \eqref{stringsOnly}, we do not consider the case where $\xi$ and $\eta$ are both even. Then, the elements of $\fg=\fpgl(m|n+2)$, where $\xi$ and $\eta$ are both odd, are represented by the vector fields
\be\label{slINk}
\begin{array}{l}
\fg_{-2}:=\Cee \xi\partial_\eta;\\
\fg_{-1}:=\Span\{\xi\partial_{x_i}  \text{~~and~~} x_i\partial_\eta \};\\
\fg_0:=\Span\{x_i\partial_{x_j} ,  \xi\partial_\xi ;  \eta\partial_\eta\}  \text{~~modulo center}.\\
\end{array}
\ee

By convention \eqref{stringsOnly}, we also have $n=0$, so in terms of supermatrices, we consider the following $\Zee$-grading of $\fh:=\fgl(m|2)$: the main block-diagonal elements are of degree 0, the block-diagonal elements just above (resp. below) the main diagonal are of degree 1 (resp. $-1$), etc.
\be\label{mat}
\begin{pmatrix}
a& r& b\\
s& A& z\\
c& u& d
\end{pmatrix}, \text{where $A\in\fgl(m)$, $a,b,c, d\in\Cee$}.
\ee

Since $\fsl(1|2)\simeq\fosp(2|2)$, the case where $n=0$, $m=1$ is already considered in Theorem~\ref{k12}.

In terms of $K_f\in\fk(1|2m)$, where $f\in\Cee[t;\xi_j, \eta_j]_{j=1, 2,\dots, m}$, see eq.~\eqref{2.2.6}, we should consider the following cases.

We realize the elements \eqref{mat} of $\fh_{\geq 0}$ by means of contact vector fields where $F:=\sum \xi_i \eta_i$
\[
\begin{array}{l}
\fh_{-2}=\text{$c\mapsto cK_1$}; \\
\fh_{-1}=\text{$s_i\mapsto s_iK_{\xi_i}$, $u_j\mapsto u_jK_{\eta_j}$};\\
\fh_0=\Span(K_t, K_{\xi_i \eta_j})_{i,j=1,2, \dots, m};\\
\fh_1=\Span\{ \eta _i(t+F),\ \xi _i(t-F)\}_{i=1, 2, \dots, m};\\
\fh_2=\Span\{(t+F)(t-F)\}.
\end{array}
\]

\ssbegin[New Bol's analogs of $i(\fpgl(2|2))$-invariant operators]{Theorem}[New Bol's analogs of $i(\fpgl(2|2))$-invariant operators]\label{T5}  The only non-zero non-scalar $i(\fpgl(2|2))$-invariant differential operators, see line~  $3$ in eq.~\eqref{d=2},
 are 
\[
B_k: {\cal F}_{-(k+3)} \tto {\cal F}_{k+1} \text{~~given by~~}  
\varphi \alpha_1^{-(k+3)/{2}} \mapsto \bar B_k(\varphi) \alpha_1^{(k+1)/{2}}  \text{for any $k\geq 0$}, 
\]
where \textup{(we use notation \eqref{comm} and \eqref{KD})}
\be
\bar B_k(\varphi)=  [D_{\eta_1}, D_{\xi_1}]_+ \circ [D_{\eta_2}, D_{\xi_2}]_+ \circ \partial_t^k(\varphi) +\nfrac{4}{(k+2)} \partial_t^{k+2}(\varphi),
 \ee
and one more operator:
\[
B: {\cal F}_{-2} \tto {\cal F}_{0}:\ \  
\varphi \alpha_1^{-1} \mapsto \bar B(\varphi),  
\] 
where
\be\label{onemore}
\bar B(\varphi)=   [D_{\eta_1},  D_{\xi_1} ]_{+}  (\varphi) + [D_{\eta_2},  D_{\xi_2}]_{+} (\varphi).
 \ee
\end{Theorem}

\begin{proof}
The proof follows the lines of the proof of Theorem \ref{k12}.  First, observe that since $i(\fpgl(2|2))$ is a~ Lie subsuperalgebra of $i(\fosp(4|2;0))$, then the operator $\bar B$ must be invariant.  
We omit the details --- a~ matter of a~ direct computation. \end{proof}

\ssbegin[New Bol's analogs of $i(\fpgl(m|2))$-invariant operators for $m\geq 3$]{Theorem}[New Bol's analogs of $i(\fpgl(m|2))$-invariant operators  for $m\geq 3$]\label{T6} There is only one non-zero non-scalar $i(\fpgl(m|2))$-invariant operator  \textup{(we use notation \eqref{comm} and \eqref{KD})}:
\[
B: {\cal F}_{-m} \tto {\cal F}_{2-m}:\ \ 
\varphi \alpha_1^{-m/2} \mapsto \bar B(\varphi)\alpha_1^{(2-m)/2}  
\] 
where  
\[
\bar B(\varphi)=  \sum_{1\leq i \leq m} [D_{\eta_i},  D_{\xi_i} ]_{+}  (\varphi) .
 \] 
\end{Theorem}

\begin{proof}
Similar to the one given in Theorem \ref{k12}.
\end{proof}

\section{The $\fpe(n+1)\subset \fm(n)$ cases} \label{BolsM} 

Since the Lie superalgebras $\fspe(n)$ are simple only for $n\geq 3$, we consider the case of subalgebra $\fpe(2)$ separately. Besides, the case of $\fm(1)$ is exceptional, see \eqref{iso}. Moreover, by convention \eqref{stringsOnly}, we consider only he case $n=1$. In the next Theorem, do not confuse the order of the operator $\partial_\tau$ which is equal to 1 with its degree which is equal to $-2$. Recall expressions \eqref{MD}.

\ssbegin[New Bol's analogs for $i(\fpe(2))$-invariant operators]{Theorem}[New Bol's analogs for $i(\fpe(2))$-invariant operators]\label{bolPe}  The only non-zero non-scalar $i(\mathfrak{pe}(2))$-invariant
linear differential operators between spaces of weighted densities are as follows.
For every $k\in\Nee$, there exists an order $k$ odd $i(\mathfrak{pe}(2))$-invariant
linear differential operator \textup{(we use notation \eqref{MD})}
\[
\begin{array}{l}B_k: {\cal F}_{\frac{k-1}{2}, \ \nu} \tto {\cal F}_{\frac{1-k}{2}, \ \nu -1+\frac{k}{2}}\text{, where} \\ 
\bar B_k =D_q(-D_\xi)^{k-1}+(k-1) (-D_\xi)^{k-2}\partial_\tau.
\end{array}
\]
\end{Theorem}

\begin{proof} Same routine as above. 
\end{proof}

\medskip

\noindent \textbf{Acknowledgments}. We heartily thank the referee for help.
S.B. and D.L. were supported by the grant NYUAD 065.

\def\eightit{\it}
\def\bib{\bf}
\bibliographystyle{amsalpha}

\begin{thebibliography}{CMZ97}

\bibitem[BaL]{BaL}
I. Basdouri  and  I. Laraiedh, The linear $\fosp(n|2)$-invariant differential operators
and cohomology, Beitr. Algebra Geom. (2014) 55: 637--643.

\bibitem[BeM]{BeM}
 M. Ben Ammar  and   W. Mtaouaa, Deformation of $\fosp(2|2)$-modules of symbols, \texttt{arXiv:1611.08147}


\bibitem[BLO]{BLO}
 N. Ben Fraj,  I. Laraiedh  and   S. Omri, Supertransvectants, cohomology, and deformations, J. of Math. Phys. 54, 023501 (2013).

\bibitem[BLO2]{BLO2}
 N. Ben Fraj,  M. Abdaoui  and   H. Raouafi, On $\fosp(2|2)$-relative cohomology of the Lie superalgebra of
contact vector fields and deformations, J. of Geometry and Physics 125 (2018) 33--48.

\bibitem[Ber]{Ber}
 J. Bernstein, Lectures on Lie Algebras. In:
\textit{Representation Theory, Complex Analysis and Integral Geometry}, Birkh\"auser (2012), 97--133


\bibitem[BL]{BL}
 J. Bernstein  and   D. Leites, Invariant differential operators and
irreducible representations of Lie superalgebras of vector fields.
Sel. Math. Sov., v.~1, no.~2, (1981), 143--160


\bibitem[Bol]{Bol}
 G. Bol, Invarianten linearer differentialgleichungen, Abh. Math. Sem. Univ. Hamburg 16 (1949) 1--28

\bibitem[BoLe]{BoLe}
F. Boniver  and  P.B.A. Lecomte, A remark about the Lie algebra of infinitesimal conformal transformations of the Euclidean space, Bull. London Math. Soc. 32:3 (2000), 263--266

\bibitem[Bc]{Bc} 
 S. Bouarroudj, Cohomology of the vector fields Lie algebras on $\mathbb{RP}^1$ acting on bilinear differential operators. Int. J. Geom. Methods Mod. Phys. 2 (2005), no. 1, 23--40. 

\bibitem[Bter]{Bter}
 S. Bouarroudj, The ternary invariant differential operators acting on the spaces of weighted densities. Theor. Math. Phys., 158, no. 2, (2009), 137--150; \texttt{arXiv:0806.1044}


\bibitem[B06]{B06}
S. Bouarroudj, Projective and conformal Schwarzian derivatives
and cohomology of Lie algebras vector of fields related to
differential operators, Int. Jour. Geom. Methods. Mod. Phys. {\bf 3} (2006), 667--696.

\bibitem[BGL]{BGL}
S.~Bouarroudj, P. Grozman  and  D. Leites, Classification of finite dimensional modular Lie superalgebras
with indecomposable Cartan matrix.
Symmetry, Integrability and Geometry: Methods and Applications
(SIGMA), 5 (2009), 060, 63 pages; \texttt{arXiv:math.RT/0710.5149}

\bibitem[BGLLS]{BGLLS}
S.~Bouarroudj, P. Grozman, A. Lebedev  and  D. Leites, I. Shchepochkina, Simple vectorial Lie algebras in characteristic
$2$ and their superizations. Symmetry, Integrability and Geometry: Methods and Applications (SIGMA) 16 (2020), 089, 101 pages; \texttt{arXiv:1510.07255}

\bibitem[BouL]{BouL}
S.~Bouarroudj  and  D. Leites, Invariant differential operators in positive characteristic. J.~Algebra. V. 499, (2018), 281--297;
\texttt{arXiv:1605.09500}

\bibitem[BLS]{BLS}
S.~Bouarroudj, D. Leites  and  I. Shchepochkina, Analogs of Bol operators for $\fpgl(a+1\vert b)\subset \fvect(a\vert b)$; \texttt{arXiv: } 


\bibitem[BO]{BO}
S.~Bouarroudj  and  V.~Ovsienko, Three cocycles on $\text{Diff}(S^1)$ generalizing the Schwarzian derivative, Int. Math. Res. Not. IMRN {\bf 1} (1998) 25--39.

\bibitem[CCLL]{CCLL}
D. Chapovalov, M. Chapovalov, A. Lebedev  and  D. Leites, The
classification of almost affine (hyperbolic) Lie superalgebras.
v.~17, Special issue in memory of F.~Berezin, (2010) 103--161;
\texttt{arXiv:0906.1860}

\bibitem[CMZ]{CMZ}
 P. B. Cohen,  Yu. Manin  and   D. Zagier, Automorphic pseudodifferential operators. Algebraic aspects of integrable systems. \textit{Progr. Nonlinear Differential Equations Appl.}, 26, Birkh\"auser Boston, Boston, MA, (1997) 17--47 

\bibitem[DGO]{DGO} 
 C. Duval, A. M. El Gradechi  and    V. Ovsienko,
Projectively and conformally invariant star-products. Comm. Math. Phys. 244 (2004), no. 1, 3--27; \texttt{arXiv:math/0301052}


\bibitem[FF]{FF}
 B. L. Feigin  and  D. B. Fuks, Invariant skew-symmetric differential operators on the line and Verma modules over the Virasoro algebra. Funct. Anal. Appl., 16:2 (1982), 114--126


\bibitem[G]{G}
 F. Gires, Conformally covariant operators on Riemann surfaces (with application to conformal and integrable models). Int. J. Modern Phys. A ({\bf 8}) (1993): 1--58.

\bibitem[GT]{GT}
 F. Gires  and  S. Theisen, Superconformally covariant operators and super W-algebras. J. Math. Phys. {\bf 34} (1993): 5964--85.



\bibitem[Gr1]{Gr1}
 P. Grozman, Classification of bilinear invariant operators on tensor
fields. Functional Anal. Appl. 14 (1980), no.~2, 127--128; for
proofs, see \texttt{arXiv:math/0509562}


\bibitem[Gr2]{Gr2}
 P. Grozman, On bilinear invariant differential operators acting on
tensor fields on the symplectic manifold. J. Nonlinear Math. Phys. 8
(2001), no. 1, 31--37; \texttt{arXiv:math/0101266}

\bibitem[GLS]{GLS}
 P. Grozman, D. Leites  and   I. Shchepochkina, Lie superalgebras of string theories.  Acta Mathematica Vietnamica,  (2001) v. 26, no.~1,
27--63, \texttt{arXiv:hep-th/9702120}

\bibitem[GLS1]{GLS1}
P. Grozman, D. Leites  and   I. Shchepochkina,  Invariant operators
on supermanifolds and standard models. Multiple facets of
quantization and supersymmetry, World Sci. Publ., River
Edge, NJ, (2002), 508--555;
\texttt{arXiv:math/0202193}.

\bibitem[H]{H}
W.-J. Huang, Superconformal covariantization of pseudodifferential operator on $(1|1)$ superspace and classical $N=2$\, $W$-superalgebras. J. Math. Phys. {\bf 35}, no. 5, (1994): 2570--82.

\bibitem[Kapp]{Kapp}
 I. Kaplansky, Graded Lie algebras I, II, preprints, Univ. Chicago,
Chicago, Ill., 1975, see
\url{http://www1.osu.cz/~zusmanovich/links/files/kaplansky/}

\bibitem[Kir]{Kir}
A.A. Kirillov, Jr. An introduction to Lie groups and Lie algebras. Cambridge Studies in Advanced Mathematics, 113. Cambridge University Press, Cambridge, 2008. xii+222 pp.

\bibitem[Kinv]{Kinv}
A.A. Kirillov, Invariant operators over geometric quantities, J. Soviet Math., 18:1 (1982), 1--21

\bibitem[IM]{IM}
K. Iohara  and  O. Mathieu, A global version of Grozman's theorem, Math. Z. 274, (2013), 955--992

\bibitem[KLV]{KLV}
Yu. Kochetkov, D. Leites  and  
A. Vaintrob, In: S.~Andima et. al. (eds.). New invariant differential operators and pseudo-(co)homology of
supermanifolds and Lie superalgebras. \textit{General Topology
and its Appl., June 1989}, Marcel Dekker, NY, 1991, 217--238

\bibitem[LO]{LO}
P.B.A. Lecomte and V. Ovsienko, Projectively invariant symbol calculus, Lett. Math. Phys. 49:3 (1999), 173--196

\bibitem[Lr]{Lr}
D. Leites, Lie superalgebras. In: Current Problems of Mathematics. Recent
developments, v. 25, 1984, VINITI, Moscow, 3--49. English
translation: J.~Soviet Mathematics, \textbf{30} (6), (1985),
2481--2512

\bibitem[LSoS]{LSoS}
D. Leites (ed.) \textit{Seminar on supersymmetry v. $1$. Algebra and
Calculus: Main chapters}, (J.~Bernstein, D.~Leites, V.~Molotkov and
V.~Shander), MCCME, Moscow, 2012, 410 pp (in Russian; a~version in
English is in preparation but available for perusal)




\bibitem[LS]{LS}
D. Leites and I. Shchepochkina, The classification of simple Lie superalgebras of vector fields. Preprint
MPIM-2003-28 \url{http://www.mpim-bonn.mpg.de/preblob/2178} (For a~short version, see 
in:  W. Nahm and L. Chau  (eds.) Differential geometric methods in theoretical physics (Davis, CA, 1988), NATO Adv. Sci. Inst. Ser. B Phys., 245, Plenum, New York, (1990), 633--651.)

\bibitem[LSh]{LSh}
D. Leites and I. Shchepochkina, Maximal graded subalgebras of simple vectorial Lie superalgebras with polynomial coefficients; \texttt{arXiv:??} 

\bibitem[MaG]{MaG}
Yu. Manin, \textit{Gauge Field Theory and Complex Geometry}. Second
edition. Springer-Verlag, Berlin, 1997. xii+346 pp.

\bibitem[N]{N}
R. Nakahama,  
Construction of intertwining operators between holomorphic discrete series representations.  
Symmetry Integrability Geom. Methods Appl. (SIGMA ) 15 (2019), Paper No. 036, 101 pp.

\bibitem[Om]{Om}
S.Omri, $\fosp(2|2)$-trivial deformations of modules of weighted
densities on the superspace $R^{1|2}$, J. Pseudo-Differ. Oper. Appl. (2015) 6:461--485.


\bibitem[OT]{OT}
V. Ovsienko and S. Tabachnikov,  \textup{
Projective differential geometry old and new. 
From the Schwarzian derivative to the cohomology of diffeomorphism groups}. Cambridge Tracts in Mathematics, 165. Cambridge University Press, Cambridge, 2005. xii+249 pp.


\bibitem[R1]{R1}
A.N. Rudakov, Irreducible representations of infinite-dimensional
Lie algebras of the Cartan type, Mathematics of the USSR-Izvestia,
v.~8 (1974), 835--866.

\bibitem[R2]{R2}
A.N. Rudakov, Irreducible representations of infinite-dimensional
Lie algebras of type $S$ and $H$, Mathematics of the USSR-Izvestia,
v.~9 (1975), 465--480.

\bibitem[Shch]{Shch}
I. Shchepochkina, How to realize Lie algebras by vector fields. Theor. Mat. Fiz., v. 147 (2006) no. 3, 821--838; \texttt{arXiv:math.RT/0509472}



\bibitem[V]{V}
O. Veblen, Differential invariants and geometry. Atti del Congr,, Int. Mat., Bologna, (1928) 6, 181--189

\end{thebibliography}

\end{document}